\documentclass{amsart}
\usepackage{graphicx}
 \usepackage[usenames,dvipsnames]{pstricks}
 \usepackage{epsfig}
 \usepackage{pst-grad} 
 \usepackage{pst-plot} 
\newtheorem{thm}{Theorem}[section]

\newtheorem{cor}[thm]{Corollary}
\newtheorem{lem}[thm]{Lemma}
\newtheorem{prop}[thm]{Proposition}
\theoremstyle{definition}
\newtheorem{defn}[thm]{Definition}
\theoremstyle{remark}

\numberwithin{equation}{section}

\renewcommand{\>}{\rangle}

\begin{document}
\title[]{Superstable  groups   acting on trees}%

\thanks{This research was supported by the "ANR" Grant  in the framework of the project "GGL"(ANR-05-JC05-47038).}%

\author{Abderezak OULD HOUCINE}%
\address{Universit\'e de Lyon;
Universit\'e Lyon 1;
INSA de Lyon, F-69621;
Ecole Centrale de Lyon;
CNRS, UMR5208, Institut Camille Jordan,
43 blvd du 11 novembre 1918,
F-69622 Villeurbanne-Cedex, France}%
\begin{abstract}   
We study superstable groups acting on trees. We prove that an action of   an $\omega$-stable group  on a simplicial tree is trivial. This shows  that an HNN-extension or a nontrivial free product with amalgamation is not $\omega$-stable.  It is also shown   that if $G$ is a superstable   group acting nontrivially on a $\Lambda$-tree, where $\Lambda=\mathbb Z$ or $\Lambda=\mathbb R$,  and if $G$  is either $\alpha$-connected and  $\Lambda=\mathbb Z$,   or if the action  is irreducible, then $G$ interprets a simple   group having a nontrivial action  on a $\Lambda$-tree.  In particular if $G$ is superstable and splits as $G=G_1*_AG_2$,  with the index of $A$  in $G_1$ different from $2$,  then $G$ interprets a simple superstable non $\omega$-stable  group. 

We will deal with "minimal" superstable groups of finite Lascar rank acting nontrivially   on $\Lambda$-trees,  where $\Lambda=\mathbb Z$ or $\Lambda=\mathbb R$.  We show that such groups $G$  have definable subgroups $H_1 \lhd H_2 \lhd G$, $H_2$ is of finite index in $G$, such that if $H_1$ is not nilpotent-by-finite then any action  of $H_1$ on a $\Lambda$-tree is trivial, and $H_2/H_1$ is either soluble or simple and acts nontrivially  on a $\Lambda$-tree.  We are interested particularly in the case where $H_2/H_1$ is simple and we show that $H_2/H_1$ has some properties similar to those of bad groups.

\end{abstract}
\email{ould@math.univ-lyon1.fr}%
\maketitle
\section{introduction}

A natural question airising  in the study of the  model theory of groups acting on trees is  to know  the model-theoretic structure of such groups, more particulary the definable subsets,  and  whether such groups can be 
stable, superstable or $\omega$-stable.  Such a study is suggested by the outstanding work of Sela \cite{Sela-Dio-V1, Sela-Dio-V2,  Sela-Dio-VI} on the elementary theory of free groups. 
It is also expected to find in this context,  new groups who have  interesting model-theoretic
 properties.  More specifically,  to find 
groups which look like  the 
so-called \textit{bad} groups.  These groups  can be seen as  the extreme counterexamples to the Cherlin-Zil$'$ber conjecture.  Roughly speaking, the conjecture and its generalisations \cite{Berline}, states that an infinite  simple group with a \textit{good} model-theoretic notion of dimension is an algebraic group over an algebraically closed field.  In the attempts to prove the conjecture, involutions play a central role. They are used, and also ideas from the classification of finite simple groups as in  Borovik's program, for proving the conjecture.  Actually, a simple group of finite Morley rank without involutions will be a counterexample. Currently, the conjecture has not been proved  even in cases where the group under consideration is linear. For more details  we refer the reader to \cite{altinel-borovik-cherlin}. 

From certain points 
of view, 
nonabelian free groups (which  act freely  on simplicial trees) have several 
properties comparable to those of bad groups.  For instance, a result announced by Bestvina and Feighn states that a proper definable subgroup of a free group is abelian, which is a property satisfied by \textit{minimal} bad groups.  Other analogies  are pointed out in \cite{Pillay}, where Pillay shows the "non CM-triviality" for free groups, a property satisfied also by infinite  simple groups of finite Morley rank \cite{Pillay-cm}.    

However,  free groups  are very far from 
being superstable. Following  \cite[Page 179]{poizat-stable}, this is proved by Gibone and we can find a proof  in \cite{poizat-type-generic}.  More generally, 
any nonabelian group satisfying the universal theory of nonabelian free groups
is not superstable \cite{Mustafin-Poizat, Ould-csa-super}, and a superstable torsion-free  hyperbolic group is cyclic \cite{Ould-csa-super}. However, Sela has shown that free groups and torsion-free hyperbolic groups are stable \cite{Sela-stab}. 

Let $G$ be a group acting on a (simplicial or real) tree. We say that the action is \emph{trivial}, if $G$ is without hyperbolic elements.

\begin{thm} \label{thm1} An action  of an  $\omega$-stable group on a simplicial  tree is trivial. 
\end{thm}

The proof of the above theorem is  simple and uses elementary facts about  actions on trees and $\omega$-stable groups. In fact,  the result  remains valid  for a large class of groups: groups whose definable abelian subgroups can be written as $B\oplus D$,  where $B$ is of finite exponent and $D$ is divisible, as  in Macintyre's theorem on the structure of $\omega$-stable abelian groups (Theorem \ref{macintyre-abelian} below).  It is worth pointing out that the precedent theorem is no longer true if we consider actions on real trees.  For instance,  $\mathbb Q$ is  of finite Morley rank and  acts freely on a real tree. 

Recall that a  group $G$ is said to split over a subgroup $A$ if it admits a nontrivial 
decomposition as an amalgamated free product $G = G_1*_AG_2$, 
 or as an 
HNN-extension $G = K*_A$.  It follows from Bass-Serre theory, that  $G$ splits over a subgroup $A$, if and only if,  $G$ acts nontrivially and  without inversions on some  simplicial tree. Hence we get the following corollary.

\begin{cor} \label{cor1}If $G$ splits over $A$ then $G$  is not  $\omega$-stable.  \qed
\end{cor}

Serre has introduced 
the notion of a group with the \textit{property FA} \cite{serre-tree}.  A group  $G$  has property FA, if $G$ has a global fixed point  whenever  $G$ acts without inversions on some simplicial  tree $T$, that is $gv=v$ for some vertex $v$ of $T$ and for every $g \in G$. It follows from Theorem \ref{thm1} and \cite[no. 6.5,  Corollary 3]{serre-tree}, that a finitely generated $\omega$-stable group has property FA. We notice that it is not known whether or not infinite finitely generated $\omega$-stable groups exist, and if they exist then this gives a counterexample to the Cherlin-Zil$'$ber conjecture in the context of $\omega$-stable groups \cite[Proposition 3.2]{ould-super-resid}.

\smallskip
A  basic example of a superstable group acting nontrivially  on a simplicial tree, we have $\mathbb Z$, which is moreover of finite Lascar rank.  More generally,  if $G$ is superstable, then the group $K=G\oplus \mathbb Z$ is also superstable, and clearly $K$ can be written as an HNN-extension, and hence acts  nontrivially on a simplicial tree.   Furthermore, if we take $G$ having a finite Lascar rank,  we get $K$ with finite Lascar rank. 

Recall that a  superstable group is \textit{$\alpha$-connected} if $G$ is connected and  $U(G)= \omega^\alpha n$, where
$\alpha$  is on ordinal, $n \in \mathbb N^*$,  $U(G)$ denotes the Lascar rank.  Recall also that actions on trees are classified  into three types:  \textit{abelian}, \textit{dihedral} and \textit{irreducible}. 

\begin{thm} \label{thm2} Let  $G$ be a superstable group acting nontrivially   on a $\Lambda$-tree, where $\Lambda=\mathbb Z$ or $\Lambda=\mathbb R$.  If $G$ is $\alpha$-connected and $\Lambda=\mathbb Z$, or if  the action is irreducible, then $G$ interprets a simple group having a nontrivial action  on a $\Lambda$-tree. \end{thm}

Combining Theorem \ref{thm1} and Theorem \ref{thm2}, we get the following. 

\begin{cor} \label{cor2} An $\alpha$-connected superstable group acting nontrivially on a simplicial tree interprets a simple superstable non $\omega$-stable  group.  \qed
\end{cor}

\begin{cor}  \label{cor3} If $G$ is superstable and splits as $G=G_1*_AG_2$, with the index of $A$  in $G_1$ different from $2$, then $G$ interprets a simple  superstable non $\omega$-stable  group acting nontrivially  on a simplicial tree.  
\end{cor}

Poizat \cite{poizat-type-generic} has shown that a nontrivial free product $G_1*G_2$ is superstable if and only if $G_1=G_2=\mathbb Z_2$.  Hence in the previous corollary, $A$ is necessarily nontrivial.

The existence of a group satisfying  conditions of Corollary \ref{cor3}, will gives a counterexample to the Cherlin-Zil$'$ber conjecture in the superstable case and to \cite[Conjecture 3 and 4]{Berline}.  Corollary \ref{cor3} is also interesting from combinatorial group theory viewpoint, as the existence of such groups will gives new examples of simple groups which can split as free products with amalgamations.  We notice that an example of a such simple group is constructed by Burger and Mozes in \cite{Burger-Mozes}, answering an old question of P. Neumann.  We also notice  that one can find superstable groups that are written as 
amalgamated free products,  but which do not satisfy  conditions of the precedent corollary. For instance, if $H$ is superstable, then $H \oplus (\mathbb Z_2*\mathbb Z_2)$ is superstable and it can be written as $(H\oplus \mathbb Z_2)*_H(H\oplus \mathbb Z_2 )$.

In the next theorem we are concerned with "minimal" superstable groups of finite Lascar rank acting on simplicial trees.  This is motivated by bad groups and by the following. Let $L$ be a non nilpotent-by-finite supertsbale group of finite Lascar rank, acting nontrivially  on a simplicial tree, and let $G \leq L$ be a definable non nilpotent-by-finite   subgroup of minimal Lascar rank  having a nontrivial action  on a simplicial  tree. Then $G$ satisfies conditions of the  theorem below.  To state the theorem we need some definitions.
  
A subgroup $H \leq G$ is said \emph{equationally-definable},  if $H$
is definable by a finite collection of  equations, that is, if
there exist  words $w_1(x), \cdots, w_n(x)$, with parameters from
$G$,  such that $H=\bigcap_{1\leq i \leq n }\{g \in G\mid G \models
w_i(g)=1\}$.  

Let $G$ be a group and $A \leq G$ be a subgroup.  The  \textit{biindex} of $A$ in  $G$ is defined to be the cardinal of the double coset space $A/G\backslash A$. Recall that a definable subset  $B$ of a stable group $G$ is said \textit{generic} if a finite number of right-cosets (or left-cosets) of $B$ covers $G$, and it is said  \textit{generous} if  $\bigcup_{g \in G}B^g$ is generic.

\begin{defn} Let $G$ be a group and $\mathcal B$ be a family of definable subgroups of $G$. We say that $\mathcal B$ is a \textit{Borel} family,  if for any $B \in \mathcal B$,  $N_G(B)/B$ is finite and $B$ is generous,  for any $g \in G$, $B^g \in \mathcal B$,   and any two elements of $\mathcal B$ are conjugate to each other. 
\end{defn}

\begin{thm} \label{thm3} Let $G$ be a  superstable group of  finite Lascar rank acting nontrivially on a $\Lambda$-tree where $\Lambda=\mathbb Z$ or $\Lambda=\mathbb R$.   Suppose that,  if $H$ is a  definable  subgroup such that $U(H)<U(G)$, and having a nontrivial action  on a $\Lambda$-tree, then $H$ is  nilpotent-by-finite.   Then  there are definable subgroups $H_1 \lhd H_2 \lhd G$ such that   $H_2$ is of finite index in $G$, and one of the following cases holds:  

\smallskip
 $(1)$   $H_1$ is connected, any action of $H_1$ on a $\Lambda$-tree is trivial,  $H_2/H_1$ is  soluble and  has a nontrivial action on a $\Lambda$-tree.

 $(2)$ $H_{2}/H_{1}$ is simple and acts nontrivially on a $\Lambda$-tree,     $H_2/H_1$ has a Borel family of  equationally-definable nilpotent subgoups  such that  there exists  $ m \in \mathbb N$ such that  for every  hyperbolic element $g$ in $H_2/H_1$, there is $1 \leq n \leq m$,  such that $g^n$ is in some $B \in \mathcal B$.  If $\Lambda=\mathbb Z$ then  $H_2/H_1=G_1*_AG_2$ with the biindex of $A$  is $2$ in both $G_1$ and $G_2$.
 \end{thm}

We will freely use notions 
of model theory and stability theory.  To the reader unfamiliar with model theory we recommend    \cite{Marker, Hodges(book)93, Poizat-cours} and for stable and superstable groups we recommend  \cite{Berline-Lascar,  poizat-stable, Pillay-book, wagner-stable}. Concerning actions on trees, we refer the reader to \cite{serre-tree, chiswell}.

The paper is organized as follows. In the next section we record the 
material that we require around actions on trees and superstable groups. 
In Section 3 we prove some lemmas on normal subgroups of groups acting on trees.   Section 4 deals with a property of elliptic elements and Section 5 is devoted to some properties of normal subgroups of superstable groups.  Section 
6 concerns  preliminary propositions and Borel famillies.  Section 7  is devoted to prove  the main results, 
and Section 8 concludes with some  remarks.

\section{Background}

We recall briefly some material
concerning  actions on trees and superstable groups. For actions on trees, we follow Serre's and Chiswell's traitement \cite{serre-tree, chiswell}.  For superstable groups,  most results can be found in \cite{Berline-Lascar}.  

\subsection{Actions on trees} $\;$

A \textit{graph} $\Gamma$ consists of a set $V$ (whose elements are called \textit{vertices}), a set $E$ (whose elements are called \textit{edges}) and a mapping
$$
E \rightarrow V \times V \times E, \quad e  \longmapsto (o(e), t(e), \bar{e}),
$$
which satisfies the following properties: for every $e \in E$, $\bar{\bar{e}}=e$, $\bar e \neq e$ and $o(e)=t(\bar  e)$.  

\smallskip
For $e \in E$, $o(e)$ is called the \textit{origin} of $e$, $t(e)$ is called the \textit{terminus} of $e$ and $\bar e$ is called the \textit{inverse}  of $e$.  The vertices $o(e)$ and $t(o)$ are often called the \textit{endpoints} of $e$.

In this context, there are natural notions of morphisms, isomorphisms, automorphisms of graphs. Let $\Gamma_1$ and $\Gamma_2$ be graphs.  A \textit{morphism} from $\Gamma_1$ to $\Gamma_2$ is a maping $f: V(\Gamma_1) \cup E(\Gamma_1) \rightarrow V(\Gamma_2) \cup E(\Gamma_2)$ such that the following properties are satisfied: 
$$
f(V(\Gamma_1)) \subseteq V(\Gamma_2),\;
 f(E(\Gamma_1)) \subseteq E(\Gamma_2), \; $$$$f(o(e))=o(f(e)), \; f(t(e))=t(f(e)), \hbox{ for any }e \in E(\Gamma_1).  
$$

With  the same manner, isomorphisms and automorphisms of graphs are defined. A \textit{simplicial tree} or a \textit{tree} is a connected graph without circuits; that is without a sequence of edges $(e_1, \cdots, e_n)$ such that $t(e_i)=o(e_{i+1})$, $o(e_i) \neq t(e_{i+1})$ and $t(e_n)=o(e_1)$.  

A group $G$ is said to acts on a tree $T$, if there is an action of $G$ by automorphisms on $T$. An element $g \in G$ is an \textit{inversion}, if $g e=\bar e$ for some edge; is \textit{elliptic} if $g v=v$ for some vertex and it is \textit{hyperbolic} if  it is neither an inversion  nor elliptic.  The action of $G$ on $T$ is said \textit{trivial} if any element $g \in G$ is either an inversion or elliptic.

Serre \cite[no. 6.1, Theorem 15]{serre-tree}  characterized countable groups with  property FA.  Bass has 
introduced the property FA$'$ for a group $G$, which is equivalent to saying that any action without inversions of $G$ on a tree is trivial \cite[Theorem 2]{alperin}. He has proved the follwing theorem which generalizes  \cite[no. 6.1, Theorem 15]{serre-tree}. 

\begin{thm} \label{thm-bass-serre} \cite{bass}\cite[Theorem 3]{alperin} Let $G$ be a group acting nontrivially and without inversions on a tree. Then either $G$ has an infinite cyclic quotient,  or $G$ splits as a  free product with amalgmation $G=G_1*_AG_2$.  \qed
\end{thm}

A \textit{real tree} is a metric space  such that any two points are joined by an unique arc  which is geodesic. We consider actions of groups by isometries on real trees.  Elliptic and hyperbolic elements are defined as above.   In that context, there is no inversions and every element is either elliptic or hyperbolic.  If a group $G$ acts without inversions on a simplicial tree $T$, then its actions on the realization of $T$ (which is a real tree), is by isometries.  For  $x$ of $T$,  we denote by $Stab(x)$ the subgroup of $G$ consisting of elements fixing $x$. 

Let $G$ acts on a  tree $T$ (real or simplicial) and $g \in G$. The \textit{hyperbolic length function} is defined by 
$$
\ell(g)=\inf \{d(x,gx) | x \in T\}. 
$$ 

Then $g$ is hyperbolic if and only if $\ell(g)>0$. 

The action is said \textit{abelian}, if  for all $g,h \in G$, $\ell(gh) \leq \ell(g)+\ell(h)$; it is said \textit{dihedral} if it is nonabelian but $\ell(gh)\leq \ell(g)+\ell(h)$, whenever $g$ and $h$ are hyperbolic; and it is said \textit{irreducible} if it is neither abelian nor dihedral.  Dihedral or abelian  actions are also often called linear actions,  because of the existence of an invariant linear subtree.

In the context of actions on $\Lambda$-trees, where $\Lambda$ is an ordered abelian group, actions on simplicial trees correspond to actions on $\mathbb Z$-trees, and actions on real trees correspond to actions on $\mathbb R$-trees.  For more details, we refer the reader to \cite{chiswell}.

The next theorem gives some equivalents caracterizations of abelian, dihedral or irreudicble actions, and it is a restatement of \cite[Proposition 3.2.7, Proposition 3.2.9, Proposition 3.3.7]{chiswell}.

\begin{thm} \label{caract-type-actions}Let $G$ be a group acting on a $\Lambda$-tree $T$. 

\smallskip
\noindent $(1)$ The following properties are equivalents: 

$(i)$ the action is abelian,

$(ii)$ the hyperbolic length function is given by $\ell(g)=|\rho(g)|$ for $g \in G$, where 

\indent $\rho : G \rightarrow \Lambda$ is  a homomorphism, 

$(iii)$ for all $g,h \in G$, $\ell([g,h])=0$. 

 \smallskip
\noindent $(2)$ The following properties are equivalents: 

$(i)$ the action is dihedral,

$(ii)$ the hyperbolic length function is given by $\ell(g)=|\rho(g)|$ for $g \in G$, where 

\indent $\rho : G \rightarrow \hbox{Isom}(\Lambda)$ is  a homomorphism whose image  contains a reflection and a 

\indent nontrivial translation, and the absolute value signs denote hyperbolic length for 

\indent the action of $\hbox{Isom}(\Lambda)$, 

$(iii)$ for all hyperbolic $g,h \in G$, $\ell([g,h])=0$, but there exists $g,h \in G$ such that 

\indent $\ell([g,h])\neq 0$. 

\smallskip
\noindent $(3)$ The following properties are equivalents: 

$(i)$ the action is irreducible,

$(ii)$ there exists hyperbolic elements $g,h \in G$ such that $\ell([g,h])\neq 0$, 

$(iii)$  $G$ contains a free subgroup of rank $2$ which acts freely, without inversions and properly discontinuously on $T$. \qed
\end{thm}

 Here, Isom$(\Lambda)$ is the group of metric automorphisms of $\Lambda$, which consists of the reflections and translations \cite[Lemma 1.2.1]{chiswell}. 

\subsection{Superstable groups.}

$\;$

\smallskip

For the definition of superstable groups and their properties, we refer the reader to \cite{Berline-Lascar}. We recall some facts that we will use. A superstable group $G$ is said to have a monomial rank if  $U(G)=\omega^\alpha n$, where $\alpha$ is an ordinal and $n \in \mathbb N$, $n \neq 0$. The following lemma was proved by Cherlin and Jaligot \cite{Cherlin-Jaligot} in the  finite Morley case and was generalized by Mustafin to superstable groups. 

\begin{lem}  \cite[Lemme 5.2]{Mustafin} \label{lem-mustafin} Let $G$ be a supertable group of monomial Lascar rank.   Let $C$ be a definable subgroup of $G$ and $X$ be a definable subset of $G$. Suppose that: 

\smallskip
$(1)$ $C$ has a finite index in its normalizer, 

$(2)$ $C \cap X$ is not generic in $C$, 

$(3)$ For every $g \in G \setminus N_G(g)$, $C \cap C^g \subseteq X$.

\smallskip
Then  $C$ is generous  in $G$.  \qed
\end{lem}

If $C$ is a definable subgroup of $G$, we let $X(C)= \bigcup_{g \in G\setminus N_G(C)} C\cap C^g$. We notice that $X(C) \leq C$,  it is definable and it is the smallest definable set satisfying condition (3) of the previous lemma. 

\begin{cor} \label{cor-generic}Let $G$ be a supertable group of monomial Lascar rank and $C$ be a definable subgroup of $G$  such that $C$ has a finite index in its normalizer. If $X(C)$ is not generic in $C$, then $C$ is generous in $G$. \qed
\end{cor}

\smallskip
The following lemma is  a particular case of \cite[Lemma 1.1.9]{wagner-stable}.  

\begin{lem} \label{lem-nilpotent-subgroups}  Let $G$ be a nilpotent superstable group. If $H<G$ is a definable subgroup of infinite index, then $N_G(H)/H$ is infinite.  \qed
\end{lem}

Recall that a group is said \textit{substable} if it is a subgroup of a stable group. The next lemma  is a sligth refinement of  \cite[Theorem 1.1.10(1)]{wagner-stable} and the proof proceeds in a similar way to the proof of that theorem.    The point which allows to make the necessary changes is that if $A$ is a subgroup then  $Z(C_G(A))$ is equationally-definable by a formula independ of $A$,  and the fact that the preimage of an equationally-definable subgroup in a quotient of $G$ by an equationally-definable subgroup, is equationally-definable. The proof is left to the reader.    

\begin{lem}\label{prop2}  Let $G$ be a substable group and $H$ be a soluble (resp. nilpotent) subgroup of $G$. Then $H$ lies in an equationally-definable soluble (resp. nilpotent) subgroup of $G$ of the same derived length (resp. nilpotency class), and the defining formula depends only on the derived length (resp. nilpotency class).  \qed 
\end{lem}

\begin{cor} \label{cor-maximal-soluble}  Let $G$ be a substable group and suppose that $G$ has maximal soluble (resp. nilpotent)  subgroups.   There exists a formula $\phi_k(x; \bar y)$ (resp. $\phi'_k(x; \bar y)$), which is an intersection of equations, such that  any maximal soluble (resp. nilpotent) subgroup  of $G$, of derived length $k$ (resp. of class $k$), is definable by  $\phi_k(x; \bar g)$ (resp. $\phi'_k(x; \bar y)$) for some $\bar g$ of $G$.  \qed
\end{cor}

\begin{thm}\label{theo-alter} \cite[Proposition 2.4]{ould-super-resid} Let  $G$ be
a superstable  group of Lascar rank $U(G)= \omega^\alpha n+\beta$ where
$\alpha$ and $\beta$ are ordinals, $n \in \mathbb N^*$ and $\beta
< \omega^\alpha$. Then either the $\alpha$-connected component
$\Gamma$ of $G$ contains a definable $\alpha$-connected normal
subgroup $K$ of $G$ such that $U(K) \geq \omega^\alpha$, or
$\Gamma$ lies in a normal definable nilpotent subgroup $K$ such
that $U(K) \geq \omega^\alpha n$. \qed
\end{thm}

The following theorem is due to A. Baudisch.  A new  proof of it can be given by using Theorem \ref{theo-alter} (see Section 8(2)). 
\begin{thm}\label{thm-Baudisch} \cite[Theorem 1.1]{baudisch} Let  $G$ be an infinite 
superstable  group.  Then there are definable subgroups $1=H_0 \lhd H_1 \lhd \dots \lhd H_n \lhd G$ such that $H_{i+1}/H_{i}$ is infinite and either  abelian or simple modulo a finite centre and $H_n$ is of finite index in $G$.  \qed
\end{thm}

\begin{cor} \label{cor-locally-soluble}\cite[Corollary 1.3]{baudisch} A locally soluble subgroup of a superstable group is soluble. \qed
\end{cor}

We end this section with the following theorem, needed elsewhere. 
\begin{thm}\label{macintyre-abelian} \cite[Th\'eor\`eme 6]{berline-abelian} An abelian  group $G$ is $\omega$-stable if and only if  $G=D \oplus B$,  where $D$ is divisible and $B$ of bounded exponent. \qed
\end{thm}

\section{Normal subgroups of groups acting on trees}

This section is devoted to  prove  some lemmas, which  deal mostely with  properties of normal subgroups of groups acting on trees. 

If $G$ is a locally nilpotent group acting without inversions and nontrivially on a tree, then it follows from Theorem \ref{thm-bass-serre} that $G$ has an infinite cyclic quotient. The following lemma  gives more informations for later use. 

\begin{lem} \label{lem-locally-nilpotent}Let $G$ be a locally nilpotent group acting nontrivially and without inversions on a tree $T$.  Then $G=\<B,s|B^s=B\>$, where $s$ is hyperbolic and $B=Stab(x)$ for some vertex $x$ of $T$. 
\end{lem}

\proof

Without loss of generality, we consider the action of $G$ on the realization of $T$ and we use results of \cite{chiswell} on group actions on real trees. By \cite[Ch4, Theorem 2.16]{chiswell}, the action of $G$ on $T$ is abelian and since $G$ has hyperbolic elements, the action is of type (III)(b) (see \cite[Page 134]{chiswell}). We denote by $A_g$ the set of elements fixed by $g$ if $g$ is elliptic, and the axis of $g$ if $g$ is hyperbolic. Let  $A= \bigcap_{g \in G}A_g$. Since the action is of type  (III)(b), $A \neq \emptyset$. Let $x \in A$ and $B=Stab(x)$. Therefore $g$ is elliptic if and only if   $g \in B$. 

Since the action is abelian, by Theorem \ref{caract-type-actions}(1), the hyperbolic length function $\ell$ is given by $\ell(g)=|\rho(g)|$, where $\rho$ is a homorphism from $G$ to $\mathbb Z$.  Therefore $B$ is the kernel of $\rho$ and  $G/B$ is infinite cyclic. Clearly the preimage $s$ of a generating element of $G/B$ is hyperbolic and thus $G=\<B,s|B^s=B\>$. \qed

\begin{lem} \label{lem-centre-action}If $G$ is a group with a nonabelian action without inversions on a $\Lambda$-tree, then every element of $Z(G)$ is elliptic. 
\end{lem}

\proof

We suppose that the assertion of the lemma is false and we show that the action is abelian. Let $g \in Z(G)$ be a hyperbolic element.  If $h_1, h_2 \in G$  are hyperbolic then $[h_1,h_2]$ is elliptic by \cite[Corollary 3.3.10]{chiswell}. Hence the action of $G$ is reducible by Theorem \ref{caract-type-actions}. 

Let $h \in G$. By \cite[Lemma 3.3.9]{chiswell}, $A_g \subseteq A_h$. Since the action is reducible, there is a linear invariant subtree $A$ by \cite[3.2.8]{chiswell}. By \cite[Lemma 3.2.10(a)]{chiswell}, $A \subseteq A_g$ and thus $A\subseteq A_h$ for any $h \in G$. It follows, by \cite[Lemma 3.2.10(c)]{chiswell}, that the action is abelian.\qed

\bigskip
In the next two lemmas, we consider the action of $G$ on the Bass-Serre tree $T$ associated to the  indicated  splitting of $G$.

\begin{lem} \label{normal-elliptic}  \label{locally-nilpotent-normal}  \label{cyclic} Let $G$ splits as a  free product with amalgmation $G=G_1*_AG_2$. If $H$ is a normal subgroup such that every element of $H$ is elliptic, then $H \leq A$.  
\end{lem}

\proof

As in the proof of the precedent lemma, without loss of generality, we consider the action of $G$ on the realization of $T$ and we use results of \cite{chiswell} on group actions on real trees. The action of $G$ on its Bass-Serre tree is either irreducible or dihedral. Hence, by \cite[Lemma 4.2.9(1)]{chiswell}, $H$ fixes some point of $T$. Therefore, $H$ is in some conjugate of $G_1$ or $G_2$. Thus, since $H$ is normal, $H \leq G_1$ or $H \leq G_2$. Again, since $H$ is normal,  we  find $H \leq A$. \qed

\smallskip
A group is said \textit{small} if it does not contain a free subgroup of rank 2.  The next lemma can certainly be extracted from  \cite{Tishin-normal}, but for completness we provide a proof. 

\begin{lem}  \label{normal-small} Let $G$ splits as a  free product with amalgmation $G=G_1*_AG_2$ such that the index of $A$ in $G_1$ different from $2$.  If $H$ is a normal subgroup of $G$ whose action by restrection is reducible, then $H \leq A$.  In particular a small normal subgroup of $G$ lies in $A$. 
\end{lem}  

\proof By Lemma \ref{normal-elliptic}, it is sufficient to show that every element of $H$ is elliptic. Suppose towards a conradiction that $H$ contains a hyperbolic element. We claim that $H$ contains a cyclically reduced (abbreviated c.r.) hyperbolic element $s$ such that every conjugate of $s$ is also c.r. 

Since the action of $H$ is reducible,  it  is either abelian or dihedral.  Let $B$ denotes the set of elliptic elements of $H$. 

Suppose that  the action is abelian. By Theorem \ref{caract-type-actions}(1),  $H=\<B,s|B^s=B\>$, where $s$ is hyperbolic. Since $H$ is normal and $s$ is a conjugate of a c.r. element,  we may assume that $s$ is c.r. and since it is hyperbolic, $s \not \in G_1 \cup G_2$. Since a conjugate of an elliptic element  is also elliptic, $B$ is a normal subgroup of $G$, and by Lemma \ref{normal-elliptic}, $B \leq A$. Hence, any conjugate of $s$ can be written as $as^{\pm 1}$ for some $a \in A$, and therefore any conjugate of $s$ is a c.r. element. This ends the proof of the claim in this case. 

Suppose that  the action is dihedral. By Theorem \ref{caract-type-actions}(1), the hyperbolic length function $\ell$ is given by $\ell(g)=|\rho(g)|$, where $\rho$ is a homorphism from $G$ to Isom($\mathbb Z$) whose image contains a reflection and a nontrivial translation.  Let $B'$ be the kernel of $\rho$. Then  $B' \leq B$ and $H/B'$ is the infinite dihedral group $D_{\infty}$.  Write $D_{\infty}=\<s_1,s_2|s_1^2=s_2^2=1\>$. Then the preimage $s$ of $s_1s_2$ is hyperbolic and $\<B',s\>=\<B',s|B'^s=B'\>$. It follows that $\<B',s\>$ is a normal subgroup.  Therefore using the same argument as in the previous case, we find that there exists a c.r. hyperbolic element $s$ such that  every conjugate of $s$ is also c.r. This ends the proof of the claim.

Let $s$ be as above, and  write $s=g_1 \cdots g_n$  in normal form.  Since $s$ is a c.r. element and hyperbolic, $n \geq 2$ and   $g_1$,  $g_n$ come from diffrent factors.  

Suppose first that $g_1 \in G_1$ and let $g \in G_1\setminus A$ such that $gA \neq g_1A$.  Let $h \in G_2 \setminus A$. Then $h^{-1}g^{-1}g_1 \cdots g_ng h$ is a conjugate of $s$ which is not cyclically reduced as $g_n \in G_2$ and $g^{-1}g_1 \not \in A$.

Suppose now that $g_1 \in G_2$ and let $g \in G_1\setminus A$ such that $Ag \neq Ag_n$. Let $h \in G_2 \setminus A$. Then $hgg_1 \cdots g_ng^{-1} h^{-1}$ is a conjugate of $s$ which is not cyclically reduced as $g_1 \in G_2$ and $g_ng^{-1} \not \in A$.  

So in each case, we find a contradiction.  This ends the proof  of the first assertion of the lemma. 

If a group $H$ acts without inversions on a tree  and the action is irreducible,    $H$ contains a free subgroup of rank 2 by  Theorem \ref{caract-type-actions}(3). Hence if $H$ is a small normal  subgroup of $G$, then the action of $H$ by restrection is reducible and thus $H \leq A$.  \qed

\begin{lem}  \label{lemme-dihedral} Let $G$ splits as a  free product with amalgmation $G=G_1*_AG_2$. Then the action of $G$ on its Bass-Serre tree is dihedral if and only if the index of $A$ is $2$ in both $G_1$ and $G_2$.  
\end{lem}
\proof
 If the action of $G$ is dihedral, then it is reducible and by Lemma  \ref{normal-small},  $A$ is of index $2$ in both $G_1$ and $G_2$. 
 
 Suppose that $A$ is of  index $2$ in both $G_1$ and $G_2$.  Let $g_1 \in G_1 \setminus A$ and $g_2 \in  G\setminus A$.  Let $H=\<A, g_1g_2\>$. Using calculation with normal forms, and since $A$ is  normal,  $H$ is also normal. Clearly the action of $H$ by restrection is abelian, and therefore, by \cite[Lemma 4.2.9(2)]{chiswell}, the action of $G$  is reducible. Since the action of $G$ on its Bass-Serre tree cannot be abelian, we deduce that it is dihedral.\qed

\begin{lem}\label{quotient1-actions} Let $G$ be a group acting nontrivially and without inversions on a $\Lambda$-tree, where $\Lambda$ is either $\mathbb Z$ or $\mathbb R$. Let $H$ be a normal subgroup such that $H \subseteq N$, where $N$ is the set of elliptic elements. Then $G/H$ acts without inversions on a $\Lambda$-tree such that  for any $g \in G$,   $g$ is hyperbolic in $G$ if and inly if $gH$ is hyperbolic in $G/H$.  In particular, if the action of $G$ is abelian (resp. dihedral, resp. irreducible), then the action of $G/H$   is abelian (resp. dihedral, resp. irreducible). 
\end{lem}
\proof

We use notations and results of \cite{chiswell}.  Let $T$ be the tree on which $G$ acts.  We give the proof only for the case $\Lambda=\mathbb R$; the case $\Lambda=\mathbb Z$ can be treated similarly and it is left to the reader. 

We  first consider  the case when the action of $G$ is abelian. By Theorem \ref{caract-type-actions}(1), there exists a homomorphism $\rho :G \rightarrow \Lambda$ whose kernel is $N$.

Let $\pi$ be the canonical  homomorphism from $G$ to $G/H$.  Let $L$ be the natural Lyndon length function of $\mathbb R$; that is the absolute value.  We define a length function $D$ on $G/H$ by $D(\pi(g))=L(\rho(g))$.  Then $D$ is well defined and it is a Lyndon length function on $G/H$.

By \cite[Theorem 4.6, Ch2]{chiswell}, $G/H$ has an action on a real tree $X$ such that $D=L_x$ for some $x \in X$, where $L_x$ denotes the Lyndon length function associated to the action at the basepoint $x$. 

Therefore, for any $g \in G$, $D(\pi(g^2))>D(\pi(g))$ if and only if $L(\rho(g^2))>L(\rho(g))$ if and only if $g \not \in N$. Hence $\pi(g)$ is hyperbolic if and only if $g$ is hyperbolic.

 By Theorem \ref{caract-type-actions}(1), an action of a group $G$ is abelian if and only if for all $g,h \in G$, $[g,h]$ is elliptic. Since the last property is satisfied by $G/H$, we find that the action of $G/H$ is abelian.

We now assume that the action of $G$  is either dihedral or irreducible. 

 By \cite[Lemma 2.9(1), Ch4]{chiswell}, $H$ has a global fixed point $a$. Hence $H$ has a bounded action.  By \cite[Lemma 2.11]{chiswell}, $G/H$ acts on the real tree $\widehat{T/H}$, in such a way that $L_{[a]}(gH)=\inf \{L_a(hg)|h \in H\}$.

Since $L_a(h)=0$ for any $h \in H$, it follows that for any $g \in G$, $L_a(hg)=L_a(g)$, and thus $L_{[a]}(gH)=L_a(g)$. 

Hence an element $gH$ is hyperbolic if and only if $L_{[a]}(g^2H)>L_{[a]}(gH)$ if and only if $L_a(g^2)>L_a(g)$. 

Using Theorem \ref{caract-type-actions}(2)-(3), we deduce that if the action of $G$ is dihedral (resp. irreducible) then the action of $G/H$ is dihedral (resp. irreducible). \qed

\section{Covering by elliptic elements}

In order to prove the conjugacy of    families of subgroups, we will need to know when  a group  acting on a tree can be covered by a finite number of translates of the set of elliptic elements.  This section is devoted to prove the following lemma. 

\begin{lem}\label{lem-hyper3}Let $G$ be a group acting without inversions on a $\Lambda$-tree with $\Lambda$ is archimedean. If there exists a finite set $S$ such  that $G=SN$, where $N$ is the set of elliptic elements,   then the action of $G$ is reducible.   \qed
\end{lem}

\noindent
\textbf{Remark}.  In the precedent lemma we connot deduce that the action is abelian. For instance,  if $G=G_1*_AG_2$ with $A$ has index $2$ in both $G_1$ and $G_2$, then the action of $G$ on the Bass-Serre tree is dihedral, and using calculation with normal forms we get that $G=SN$ for some finite subset $S$. 

Lemma \ref{lem-hyper3} is a consequence of the next two lemmas. 

\begin{lem} \label{lem-hyper}Let $G$ be a group acting without inversions on a $\Lambda$-tree with $\Lambda$ is archimedean. If $a,b \in G$ are hyperbolic, then there exists $p \in \mathbb N$ such that for any $n \geq p$, $ab^n$ is hyperbolic. 
\end{lem}

\proof

We use notations and results of \cite{chiswell}. We first prove the following claim.   

\smallskip
\noindent
\textit{Claim. Suppose that $\ell(b)> \ell(a)$. Then there exists $p \in \mathbb N$ such that for any $n \geq p$, $ab^n$ is hyperbolic. }
\smallskip

\proof

If $A_{a} \cap A_{b}=\emptyset$ or if $a$ and $b$ meet coherently,  then  by \cite[Lemma 3.2.2]{chiswell} or by  \cite[Lemma 3.3.1]{chiswell} we have for $n \geq 1$,  $$\ell(ab^n)\geq \ell(a)+\ell(b^n),$$ and thus $ab^n$ is hyperbolic for $n \geq 1$.

In the remainder of this proof we assume  that $a$ and $b^{-1}$ meet coherently. We  treat the following three cases. 

\smallskip
\textit{Case 1}.  $\Delta(a,b)<\min\{\ell(a), \ell(b)\}$. 

In particular we have $\Delta(a,b^n)<\min\{\ell(a), \ell(b^n)\}$, for  $n \geq 1$. Since $a$ and $b^{-n}$ meet coherently, it follows by \cite[Lemma 3.3.3(2)]{chiswell}, 
$$
\ell(ab^n)=\ell(a)+\ell(b^n)-2 \Delta>0,
$$
and thus $ab^n$ is hyperbolic for  $n \geq 1$.

\smallskip
\textit{Case 2}.  $\Delta(a,b)>\min\{\ell(a), \ell(b)\}=\ell(a)$. 

\smallskip
In particular we have $\Delta(a,b^n)>\min\{\ell(a), \ell(b^n)\}$, for  $n \geq 1$. 
Let $n \geq 1$. Since $a$ and $b^{-n}$ meet coherently,  $a^{-1}$ and $b^n$ meet coherently. Therefore, by \cite[Lemma 3.3.4(2)]{chiswell}, 
$$
\ell(b^na)=\ell(b^n)-\ell(a)>0, 
$$
and thus $ab^n$, which is a conjugate of $b^na$,  is hyperbolic for  $n \geq 1$.

\smallskip
\textit{Case 3}.  $\Delta(a,b)=\min\{\ell(a), \ell(b)\}=\ell(a)$.

\smallskip

Suppose that there exists  $m \geq 1$ such  that $b^ma$ is hyperbolic. It follows, by \cite[Lemma 3.3.5(2)(1)]{chiswell}  applied to the tuple $(b^m, a^{-1})$, that $\Delta(b^ma,b)=\ell(b^ma)>0$,  and $b^ma$ meet $b$ coherently. Hence, by \cite[Lemma 3.3.2]{chiswell}, $b^na$ is hyperbolic for any $n \geq m$. Therefore,  $ab^n$, which is a conjugate of $b^na$,  is hyperbolic for any $n \geq m$.  

Therefore,  it suffices to show that   one of the elements  $ba$, $b^2a$,  $b^3a$ is hyperbolic.  We suppose that  $ba$ and $b^2a$ are elliptic and we show that $b^3a$ is hyperbolic.  For this case, the situation is illustrated by the picture below. 

Let $r$ be  the right-hand of $A_{b} \cap A_{a^{-1}}$ and $l$ be the left-hand of  $A_{b} \cap A_{a^{-1}}$.

 By \cite[Lemma 3.3.5(3)]{chiswell} applied  to the tuples $(b, a^{-1}), (b^2, a^{-1})$, we have 
 $$A_{b^ia} \cap A_b=\{q_i\},$$
 where $q_i$ is the midpoint of $[r, b^iar]$,  for $i \in \{1,2\}$.  Since $ar=l \in A_b$, we find  $bar \neq b^2ar$. Therefore $A_{b^2a} \cap A_{ba} =\emptyset$. Further,  $[q_2,q_1]$ is the bridge between $A_{b^2a}$ and $A_{ba}$.

 By \cite[Lemma 3.2.2]{chiswell} applied to the tuple $(b^2a, ba)$,  $b^2aba$ is hyperbolic and 
 $$
 [q_2,q_1] \subseteq A_{b^2aba}. 
 $$
 
 If $A_{b^2aba} \cap A_{a^{-1}} =\emptyset$, then, by \cite[Lemma 3.2.2]{chiswell},  $b^2abaa^{-1}$ is hyperbolic and thus $ab^3$ is hyperbolic. 
 
 Suppose that $A_{b^2aba} \cap A_{a^{-1}} \neq \emptyset$. Since 
 $$
 b^2abaq_1=bq_1,
 $$
 we find that $b^2aba$ translates in the same direction as $b$. Since $b$ and $a^{-1}$ meet coherently, it follows that $b^2aba$ and $a^{-1}$ meet coherently. Hence, by \cite[Lemma 3.3.1]{chiswell} $b^2ab$ is hyperbolic and thus $ab^3$ is hyperbolic.  This completes the proof of the  claim. \qed

Since $\Lambda$ is archimedean, we can find $q \in \mathbb N$ such that $\ell(b^q)>\ell(a)$.  By applying the above claim to $(a,b^q)$ we get the required conclusion.\qed

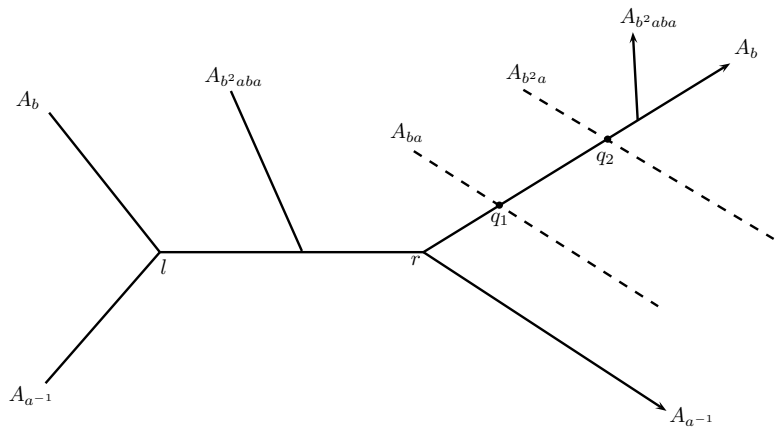
\begin{center}

\scalebox{0.8} 
{
\begin{pspicture}(0,-3.53)(13.44,3.53)
\psline[linewidth=0.04,arrowsize=0.05291667cm 2.0,arrowlength=1.4,arrowinset=0.4]{->}(1.32,1.75)(3.16,-0.57)(7.54,-0.57)(12.64,2.57)
\psline[linewidth=0.04cm](3.16,-0.57)(1.26,-2.75)
\psline[linewidth=0.04cm,arrowsize=0.05291667cm 2.0,arrowlength=1.4,arrowinset=0.4]{->}(7.54,-0.57)(11.58,-3.21)
\psline[linewidth=0.04cm,linestyle=dashed,dash=0.16cm 0.16cm](9.2,2.13)(13.36,-0.35)
\psline[linewidth=0.04cm,linestyle=dashed,dash=0.16cm 0.16cm](7.38,1.11)(11.44,-1.47)
\usefont{T1}{ptm}{m}{n}
\rput(12.91,2.815){$A_b$}
\usefont{T1}{ptm}{m}{n}
\rput(12.01,-3.305){$A_{a^{-1}}$}
\usefont{T1}{ptm}{m}{n}
\rput(7.42,-0.705){$r$}
\usefont{T1}{ptm}{m}{n}
\rput(3.23,-0.805){$l$}
\usefont{T1}{ptm}{m}{n}
\rput(0.97,1.955){$A_b$}
\usefont{T1}{ptm}{m}{n}
\rput(1.03,-2.965){$A_{a^{-1}}$}
\psdots[dotsize=0.12](8.8,0.21)
\psdots[dotsize=0.12](10.6,1.31)
\usefont{T1}{ptm}{m}{n}
\rput(8.81,-0.025){$q_1$}
\usefont{T1}{ptm}{m}{n}
\rput(10.57,0.995){$q_2$}
\usefont{T1}{ptm}{m}{n}
\rput(9.26,2.395){$A_{b^2a}$}
\usefont{T1}{ptm}{m}{n}
\rput(7.27,1.415){$A_{ba}$}
\psline[linewidth=0.04cm](5.52,-0.55)(4.34,2.11)
\psline[linewidth=0.04cm,arrowsize=0.05291667cm 2.0,arrowlength=1.4,arrowinset=0.4]{->}(11.1,1.63)(11.02,3.09)
\usefont{T1}{ptm}{m}{n}
\rput(11.29,3.335){$A_{b^2aba}$}
\usefont{T1}{ptm}{m}{n}
\rput(4.39,2.355){$A_{b^2aba}$}
\end{pspicture} 
}

\end{center}

\begin{lem} \label{lem-hyper2}Let $\Lambda$ be an archimedean ordered abelian group and  $G$ be a group acting without inversions and irreducibly on a $\Lambda$-tree. If $S \subseteq G$ is  a finite set,  then there exists a hyperbolic element $g$ such that for any $s \in S$, $sg$ is hyperbolic. 
\end{lem}

\proof

Since the action of $G$ is irreducible, we find $a,b \in G$  hyperbolic such that $A_a \cap A_b =\emptyset$.  We first  show  the following claim.

 \smallskip
\noindent
\textit{Claim.  Let $g \in G$. Then there exists $p \geq 1$ such that for any $m \geq p$, there exists $q_m$ such that for any $n \geq q_m$, $ga^mb^n$ is hyperbolic.  }
\smallskip

\proof 

We treat the following two cases. 

\smallskip
\textit{Case 1. $ga$ is hyperbolic. } 

\smallskip
By Lemma \ref{lem-hyper}, there exists $p \geq 1$ such that for any $m \geq p$, $ga^m$ is hyperbolic. Again by  Lemma \ref{lem-hyper}, for any $m \geq p$, there exists $q_m$ such that for any $n \geq q_m$, $ga^mb^n$ is hyperbolic. 

\smallskip
\textit{Case 2. $ga$ is elliptic. } 

\smallskip
If $A_{ga} \cap A_a = \emptyset$, then by \cite[Lemma 3.2.3]{chiswell}, $ga^2$ is hyperbolic. Argument similar to that in Case 1, gives the desired result. 

So we assume  that $A_{ga} \cap A_a \neq  \emptyset$. 

If $\Delta(ga, a)>0$, since $ga$ is elliptic, $ga$ and $a$ meet coherently and thus by \cite[Lemma 3.3.2]{chiswell}, $ga^2$ is hyperbolic.   The conclusion follows as in Case 1.

So we assume that $A_{ga} \cap A_a= \{r\}$. If $A_b \cap A_{ga}=\emptyset$ or $A_b \cap A_{ga}\neq \emptyset$ and $\Delta (b,ga)>0$,  then  by \cite[Lemma 3.2.2]{chiswell} or by  \cite[Lemma 3.3.1]{chiswell}, we find $gab$ hyperbolic. 
Thus $bga$ is hyperbolic. Hence, by Lemma \ref{lem-hyper}, there exists $p$ such that for any $m \geq p$, $bga^m$ is hyperbolic. Reappling the same lemma, we find that for any $m \geq p$, there exists $q_m$ such that for any $n \geq q_m$,  $ga^mb^n$ is hyperbolic.

So we assume that $A_b \cap A_{ga}=\{s\}$. 

Since $a^{-1}g^{-1}$ is elliptic, $b$ and $a^{-1}g^{-1}$ meet coherently. 

If $bga$ is hyperbolic, then $gab$ is hyperbolic and the proof runs  as above.

So we assume that $bga$ is elliptic. By \cite[Lemma 3.3.5(3)]{chiswell}, since $\Delta(b,ga)=0=\ell(ga)$ and $\ell(b)>\ell(ga)$,  we find $d(A_{bga}, A_{ga})>0$. Hence $A_{bga} \cap A_{ga}=\emptyset$.

By \cite[Lemma 3.2.2]{chiswell} applied to the tuple $(bga, (ga)^{-1})$,  

$A_{b} \cap A_{ga}=[q, (ga)^{-1}q]=\{q\},$
where $q \in A_{ga}$.  We conclude that
$$
A_{ga} \cap A_a=\{r\}, \quad A_{ga} \cap A_b=\{s\}, \quad A_{bga} \cap A_b=\{q\},  \quad A_{bga} \cap A_{ga}=\emptyset, $$
and thus 
$$
A_{bga} \cap A_a=\emptyset . 
$$

Therefore by \cite[Lemma 3.2.2]{chiswell}, $bga^2$ is hyperbolic. The conclusion follows as above. This ends the proof of the claim. \qed

Let $S=\{s_1,\dots, s_n\} \subseteq G$ be a finite set.  For each $s_i \in S$, let $p_i$ be the integer given by the above claim.  Let $p=\max \{p_i| 1 \leq i \leq n\}$.  Hence there exists $q_1(p), \dots, q_n(p)$ such that for any $n \geq q=\max\{q_i(p)| 1 \leq i \leq n\}$, $s_ia^pb^n$ is hyperbolic. Therefore,  by taking $g=a^pb^q$, which is hyperbolic,  we get the desired conclusion. \qed

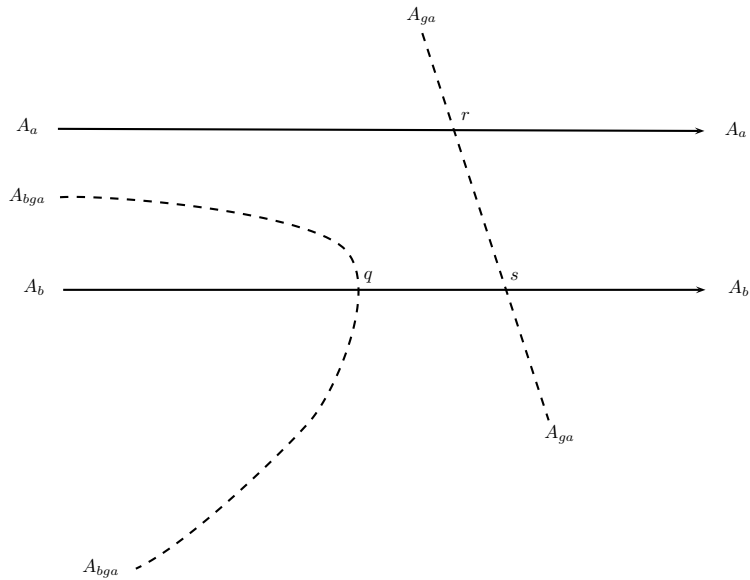
\begin{center}
\scalebox{0.7} 
{
\begin{pspicture}(0,-5.46)(14.84,5.46)
\psline[linewidth=0.04cm,arrowsize=0.05291667cm 2.0,arrowlength=1.4,arrowinset=0.4]{->}(1.38,3.12)(13.66,3.08)
\psline[linewidth=0.04cm,arrowsize=0.05291667cm 2.0,arrowlength=1.4,arrowinset=0.4]{->}(1.48,0.06)(13.68,0.06)
\psline[linewidth=0.04cm,linestyle=dashed,dash=0.16cm 0.16cm](8.3,4.94)(10.7,-2.42)
\psbezier[linewidth=0.04,linestyle=dashed,dash=0.16cm 0.16cm](1.42,1.82)(2.26,1.88)(6.217083,1.5486928)(6.86,0.8)(7.5029173,0.051307186)(6.64,-1.88)(6.18,-2.42)(5.72,-2.96)(3.84,-4.8)(2.86,-5.24)
\usefont{T1}{ptm}{m}{n}
\rput(14.26,3.105){$A_a$}
\usefont{T1}{ptm}{m}{n}
\rput(14.31,0.085){$A_b$}
\usefont{T1}{ptm}{m}{n}
\rput(10.05,0.325){$s$}
\usefont{T1}{ptm}{m}{n}
\rput(9.12,3.365){$r$}
\usefont{T1}{ptm}{m}{n}
\rput(7.27,0.305){$q$}
\usefont{T1}{ptm}{m}{n}
\rput(8.29,5.265){$A_{ga}$}
\usefont{T1}{ptm}{m}{n}
\rput(2.2,-5.235){$A_{bga}$}
\usefont{T1}{ptm}{m}{n}
\rput(10.91,-2.695){$A_{ga}$}
\usefont{T1}{ptm}{m}{n}
\rput(0.8,3.165){$A_a$}
\usefont{T1}{ptm}{m}{n}
\rput(0.93,0.105){$A_b$}
\usefont{T1}{ptm}{m}{n}
\rput(0.8,1.825){$A_{bga}$}
\end{pspicture} 

}\end{center}

\section{Normal subgroups of superstable groups}

\smallskip
We now trun to some properties of normal subgroups of superstable groups. 

\begin{lem} \label{thm-connected} Let $G$ be a superstable group of finite Lascar rank. Then $G$ has a definable connected normal subgroup $H$ such that   $G/H$ is nilpotent-by-finite. Furthermore we can take $H$ maximal among normal connected subgroups of $G$ and  $H \leq G^{\circ}$. 
\end{lem}

\proof

Let $(H_i)_{i \in \lambda}$ be the list of connected normal subgroups  of $G$.  Each $H_i$ is connected definable and $0$-indecomposable. By the Indecomposablity Theorem \cite[V, Theorem 3.1]{Berline-Lascar}, the subgroup $H$ generated by the family $(H_i)_{i \in \lambda}$ is definable and connected. Clearly $H$ is normal and maximal among   normal connected subgroups  of $G$.  We notice also that $H \leq G^\circ$.

Let $\bar G=G/H$. If $\bar G$ is finite, the result is clear. We assume that $G/H$ is infinite, and thus $U(\bar G)=\omega^\alpha n+\beta, n \in \mathbb N, n \neq 0$, $\alpha=0$, $\beta=0$. 

By Theorem \ref{theo-alter}, either the $0$-connected omponent $(\bar G)^\circ$  of $\bar G$ contains a definable $0$-connected normal
subgroup $\bar K$ of $\bar G$ such that $U(\bar K) \geq \omega ^0=1$, or
$(\bar G)^\circ$ lies in a normal definable nilpotent subgroup $\bar K$ such
that $U(\bar K) \geq \omega^\alpha n=n=U(\bar G)$ (In our contexet $0$-connectedness coincides with connectedness). 

We show that the first case is impossible. Let  $\bar K \leq (\bar G)^\circ$  connected and infinite.  Its preimage $K$ is connected definable infinite and $H \leq K$. Since $(\bar G)^\circ=G^\circ /H$, we get $H \leq K \leq G^\circ$, contradicting the maximality of $H$. 

Hence the first case is impossible and therefore $(\bar G)^\circ$ lies in a definable nilpotent subgroup having the same Lascar rank as $\bar G$, and thus $\bar G$ is nilpotent-by-finite.  \qed

\begin{cor} \label{quotient-par-la-composante-connexe}If $G$ is a superstable group of finite Lascar rank then $G/G^{\circ}$ is nilpotent-by-finite.\qed \end{cor}

\smallskip
\noindent
\textbf{Remark}.  Analogously, if $G$ is superstable then $G/G^\circ$ is soluble-by-finite (see Section 8(3)).

\begin{lem} \label{propo2}Let $G$ be an  $\alpha$-connected supserstable group. Let $N$ be a normal subgroup of $G$. Then either $N \leq Z(G)$ or $N$ contains an infinite $\alpha$-connected definable normal subgroup $H$ which is normal in $G$.  In the later case we can take $H$ maximal among normal $\alpha$-connected subgroups $L$ contained in $N$. 
\end{lem}

\proof

Assume that $N \not \leq Z(G)$ and  let $a \in N \setminus Z(G)$. If $a^G$ is finite, then $C_G(a)$ is of  finite index, and since $G$ is connected, we conclude that $a \in Z(G)$. 

Hence $a^G$ is infinite.  By \cite[VI, Lemma 2.3]{Berline-Lascar}, $a^G$ is $\alpha$-indecomposable. Therefore $a^{-1}a^G$ is also $\alpha$-indecomposable. Let $X=a^{-1}a^G$, and $H$ the subgroup generated by $\{X^g|g \in G\}$. By the Indecomposablity Theorem \cite[V, Theorem 3.1]{Berline-Lascar}, $H$ is definable and $\alpha$-connected. Then $H$ is an infinite $\alpha$-connected normal subgroup of $G$. 

Concerning the maximality of $H$, the proof proceeds in a similar way to that of  Lemma \ref{thm-connected}. \qed

\begin{lem} \label{cor-maximal-normal} If $G$ is an $\alpha$-connected superstable group then any maximal normal subgroup $N$ of $G$ is either definable or has a finite index in $G$. \end{lem}

\proof

The proof is by induction on $U(G)$. Let $N$ be a maximal normal subgroup of $G$. If $N \leq Z(G)$, then either $N=Z(G)$, in which case $N$ is definable, or $G$ is abelian. In the later case we find that $G/N\simeq \mathbb Z/ p \mathbb Z$ for some prime $p \geq 2$ and thus $N$ has a finite index.

If $N \not \leq Z(G)$, then, by Lemma \ref{propo2}, $N$ has an infinite normal $\alpha$-connected definable subgroup $H$. Hence $G/H$ is $\alpha$-connected and   $U(G/H)<U(G)$. The image of  $N$ in $G/H$ is also a maximal normal subgroup, and by indcution it is either definable or has a finite index. Hence $N$ is either  definable or has a finite index.  \qed

\begin{lem} \label{lem-finite-index}Let $G$ be an $\alpha$-connected superstable group and $N$ be a normal subgroup of finite index. Then $G$ has a definable normal subgroup $K \leq N$ such that $G/K$ is abelian. 
\end{lem}

\proof

We proceed by induction on $U(G)$.  If $N \leq Z(G)$, $G$ is abelian and the result is celar. Se we assume that $N \not \leq Z(G)$, and thus by Lemma \ref{propo2}, there exists an infinite $\alpha$-connected normal subgroup $H \leq N$. We have  $U(G/H)<U(G)$, $G/H$ is $\alpha$-connected  and the image of $N$ in $G/H$ has a finite index. The result follows by induction. \qed

\section{Preliminary propositions \& Borel famiilies}

This section is devoted to prove some preparatory propositions. Some of them are related to  irreducible actions of superstable groups. The remaining is devoted to find a Borel family in "minimal" superstable groups of finite Lascar rank,  acting on trees.

\begin{prop}  \label{propo-quotient} If $G$ is a superstable group  which has an infinite cyclic quotient or an infinite dihedral quotient, then $G$ is not $\alpha$-connected. 
\end{prop}

\proof

Suppose towards a contradiction that $G$ is $\alpha$-connected. We may assume that $G$ has a minimal Lascar rank. 

Assume  that $G$ has an infinite cyclic quotient and let $N$ be a normal subgroup such that $G/N$ is infinite cyclic. If $N \not \leq Z(G)$, by Lemma \ref{propo2}, there is a maximal definable $\alpha$-connected subgroup $H \leq N$ which is normal in $G$.   By \cite[Corollary 8.6, III. 8]{Berline-Lascar}, and \cite[Lemma 1.4, V.1]{Berline-Lascar}, $G/H$ is also $\alpha$-connected. But $G/H$ has an infinite cyclic quotient and $U(G/H)<U(G)$;  a contradiction with the fact that $G$ has a minimal Lascar rank.  
Hence $N \leq Z(G)$ and thus $G=N \oplus \mathbb Z$.  

We now show that a superstable abelian  group $K$ which can be written as $N \oplus \mathbb Z$ is not connected. Suppose towards a contradiction that $K$ is connected.  Let $K_n=n! K$. Then $(K_n| n \in \mathbb N)$ is a descending chain of definable subgroups each of which is connected, and since $K$ is superstable we get $K_{n+1}=K_n$ for some $n \in \mathbb N$. But, since $K_n=n!N\oplus n! \mathbb Z=(n+1)!N\oplus (n+1)! \mathbb Z$, we find a contradiction. This ends the proof  when $G$ has an infinite cyclic quotient. 

Assume now  that   $G$ has an infinite dihedral quotient and let $N$ be a normal subgroup such that $G/N$ is infinite dihedral. By the same methode as above,  we find $N \leq Z(G)$. Since $Z(D_{\infty})=1$, we deduce that $N=Z(G)$.  Since $G$ is connected, $G/Z(G)$ is connected.  Write $D_{\infty}=\<a|a^2=1\>*\<b|b^2=1\>$. Then $C_G(ab)$ has a finite index in $G/Z(G)$; a contradiction with  the connectedness of $G/Z(G)$.   \qed

If $G$ is acts  nontrivially on a (simplicial) tree and if the action is reducible, then   $G$ has a  cyclic or a dihedral infinite quotient by Theorem \ref{caract-type-actions}.    Hence we get the following. 

\begin{cor} \label{cor-irreducible} A nontrivial action of a superstable  $\alpha$-connected group  on a simplicial tree is irreducible. \qed
\end{cor}

The precedent corollary is not longer true if we consider actions on real trees. Indeed, $\mathbb Q$ is a connected group of finite Morely rank, and $\mathbb Q$ acts (freely) on a real tree, were the action is clearly abelian. However, the following proposition gives a reduction. 

\begin{prop} \label{reducible-actions-on-reel-trees}A reducible action of a superstable $\alpha$-connected group $G$ on a real tree is abelian.  Furthermore if the action is nontrivial then $G$ has a definable normal subgroup $K \leq N$, where $N$ is the set of elliptic elements, such that $G/K$ is metabelian; more precisely,  $N/K \leq Z(G/K)$ and $(G/K)/(N/K)$ is a subgroup of $\mathbb R$. \end{prop}
\proof
Suppose towards a contradiction that $G$ is $\alpha$-connected with a dihedral action on a real tree. We assume that $G$ has a minimal Lascar rank. 

We claim that we may assume $Z(G)=1$.  If $U(Z(G))\geq \omega^\alpha$ then $U(G/Z(G))<U(G)$. By Lemma \ref{lem-centre-action}, every element of $Z(G)$ is elliptic and by Lemma \ref{quotient1-actions}, $G/Z(G)$ has a dihedral action on a real tree; a contradiction with the minimality of the Lascar rank of $G$. 

Hence $U(Z(G))<\omega^\alpha$ and by \cite[Lemma 1.4]{Berline}, $G/Z(G)$ is centerless. By Lemma \ref{quotient1-actions} $G/Z(G)$ has a dihedral action on a real tree. Thus by replacing $G$ by $G/Z(G)$ we may assume that $G$ is centerless as claimed.

 By Theorem \ref{caract-type-actions}(2), $G$ has a normal subgroup $H$ such that $H \subseteq N$, where $N$ is the set of elliptic elements, and $G/H$ is a split extension of an infinite abelian group by $\mathbb Z_2$.
 
 If $H$ is finite then $G$ is finite-by-abelian-by-$\mathbb Z_2$, and thus, since $G$ is connected, $G$ is abelian. Hence the action of $G$ is abelian; a contradiction with our assumption. 
 
 Therefore $H$ is infinite.  Since $Z(G)=1$, by Lemma \ref{propo2} $G$ has a normal definable $\alpha$-connected  subgroup $L$ with $L \leq H$. By Lemma  \ref{quotient1-actions} $G/L$ has a dihedral action on a real tree. Since $U(G/L)<U(G)$ we get also a contradiction. 
 
 We now prove the last assertion of the proposition and we proceed by induction.  Since the action is abelian, $N$ is a normal subgroup and $G/N$ is a subgroup of $\mathbb R$ by Theorem \ref{caract-type-actions}. If $N \leq Z(G)$ then $G$ is metabelian and there is no thing to prove. Otherwise, by  Lemma \ref{propo2}, $G$ has a normal definable  $\alpha$-connected subgroup $H\leq N$. Now $U(G/H)<U(G)$ and by Lemma \ref{quotient1-actions}, $G/H$ has an abelian action on a real tree. The result follows by induction.   \qed

\begin{prop} \label{irreducibles-actions}Let $G$ be a superstable group with an irreducible action on a $\Lambda$-tree, where $\Lambda$ is either $\mathbb Z$ or $\mathbb R$.  Let $N$ be the set of elliptic elements.  Suppose that   if $H$ is a definable subnormal subgroup such that $U(H)<U(G)$,  then the action of $H$ is reducible.  Then $G$ has a serie $H_1 \lhd H_2 \lhd G$ such that $H_1 \subseteq N$, $H_2$ has a finite index in $G$ and $H_2/H_1$ is simple and acts nontrivially on a $\Lambda$-tree. 
\end{prop}

\proof 

Let $K_0=1 \lhd K_2 \lhd \dots \lhd K_n\lhd G$ be the serie given by Theorem  \ref{thm-Baudisch}.  Since $K_n$ is of finite index, its action by restrection is irreducible. Since $K_n/K_{n-1}$ is infinite, $U(K_{n-1})<U(G)$ and thus the action of $K_{n-1}$ is reducible. It follows by \cite[Lemma 4.2.9]{chiswell}, that $K_{n-1} \subseteq N$. Hence $K_{n}/K_{n-1}$ has an irreducible action on a real tree by Lemma \ref{quotient1-actions}. Therefore $K_n/K_{n-1}$ is nonabelian and thus it is simple modulo a finite center. By taking $H_1$ to be the preimage of $Z(K_n/K_{n-1})$ and $H_2=K_n$,  we get the desired conclusion. \qed

In the sequel, we seek a Borel family, but before doing so, 
we   prove the next proposition which  gives sufficient conditions for a family to be a Borel family.

\begin{defn} Let $G$ be a group and $\mathcal B$ a family of subgroups of $G$. We say that $\mathcal B$ is a \textit{ferment} family if any $B \in \mathcal B$ is definable,  $N_G(B)/B$ is finite and for any $g \in G$, $B^g \in \mathcal B$. 
\end{defn}

\begin{prop} \label{prop-conjuguaison-borel} Let $G$ be an $\alpha$-connected  superstable group acting nontrivially on a $\Lambda$-tree,   such that a finite number of right-translates of the set of elliptic elements  cannot cover $G$. Let $\mathcal B$ be a ferment family  satisfying the following conditions:

\smallskip

$(1)$ the action, by restrection,  of every $B \in \mathcal B$ is abelian and nontrivial, 

$(2)$ every $B \in \mathcal B$ contains a hyperbolic element, 

$(3)$ if $B_1, B_2 \in \mathcal B$, $B_1 \neq B_2$, then $B_1 \cap B_2$ is without hyperbolic elements. 

\smallskip
Then every $B \in \mathcal B$ is generous and any two elements of $\mathcal B$ are conjugate to each other. 

\end{prop}
\begin{defn} If $\mathcal B$ is a ferment family of $G$ which satisfies the conclusion of the above proposition, then we say that $\mathcal B$ is a \textit{Borel} family. 
\end{defn}

\proof

Let $B \in \mathcal B$. We show first that $B$ is generous. By Corollary \ref{cor-generic}, it suffices to show that $X(B)$ is not generic in $B$.  Since the action of $B$ by restrection is abelian and nontrivial, by Theorem \ref{caract-type-actions}(1), the set of elliptic elements $N$ of $B$  is a normal subgroup and $B/N$ is infinite.

Let $g \in G\setminus N_G(B)$.  By (2), $B \cap B^g \leq N$ and thus $X(B) \leq N$. Since $B/N$ is infinite, $X(B)$ is not generic in $B$. 

We now  show that if $B_1, B_2 \in \mathcal B$, such that $B_1 \neq B_2$, then $B_1$ and $B_2$ are conjugate. Let 
$$
X_i= \bigcup_{g \in G} B_i^g, \hbox{ for } i=1,2. 
$$

Since $G$ is connected and $X_1$ and $X_2$ are generic, we find that $X_1 \cap X_2$ is generic. 

Suppose towards a contradiction that  for any $(g,g') \in G^2$,  $B_1^g \cap B_2^{g'}$ does not contain a hyperbolic element.   We have 
$$
X_1 \cap X_2 =  \bigcup_{(g,g') \in G^2} B_1^g \cap B_2^{g'}, 
$$
and for any $(g,g') \in G^2$,  $B_1^g \cap B_2^{g'} \leq N$, where $N$ is the set of elliptic elements of $G$. Therefore $X_1 \cap X_2 \leq N$ and thus $G = a_1N \cup \dots \cup a_nN$; a contradiction to  the hypothesis of the proposition. 

Therefore, for some $(g,g') \in G^2$,  $B_1^g \cap B_2^{g'}$ contains a hyperbolic element and thus $B_1^g=B_2^{g'}$. This ends the proof of the proposition. \qed

\begin{cor} Let $G$ be an  $\alpha$-connected superstable group with an irreducible action  on a real tree. If $\mathcal B$ is a ferment family of $\alpha$-connected small subgroups satisfying conditions (2) and (3)  of Proposition \ref{prop-conjuguaison-borel}, then $\mathcal B$ is a Borel family. \end{cor}

\proof 

It is a consequence of Proposition  \ref{prop-conjuguaison-borel}, Proposition \ref{reducible-actions-on-reel-trees}, Lemma \ref {lem-hyper3} and of the fact that an action of a small group is reducible. \qed

\begin{defn}  Let $G$ be a group of finite Lascar rank $n$. It follows from \cite[Corollary VI. 2.12]{Berline-Lascar} that  any soluble (resp. nilpotent) subgroup is of derived length (resp. class) at most $n$.  We let,  for  $0 \leq k \leq n$,  $\phi_k(x; \bar y)$ to be the formula given by Corollary \ref{cor-maximal-soluble}, and for $n+1 \leq k \leq 2n$, $\phi_k(x; \bar y)$ to be the formula $\phi'_{k-n}(x; \bar y)$ given by Corollary \ref{cor-maximal-soluble}.  For $0 \leq k \leq 2n$, we let $\varphi_k(x; z, \bar y):=\phi_k(x^z; \bar y)$. 

Let $B$ be a soluble subgroup of $G$. We let $B_k$ to be intersection of  $\varphi_k$-definable  subgroups $K$ of $G$,  such that $K \cap B$ has a finite index in $B$. By Baldwin-Saxl's lemma, $B_k$ is a finite intersection. We let $B^*$ to be $B \cap (\cap_{0 \leq k \leq 2n}B_k)$.   We say that $B$ is $*$-connected if $B=B^*$.  By Lemma \ref{prop2},   $B^*$ exists and has a finite index in $B$.  We notice that $B^*$ is  not necessarily  definable. 

\end{defn}

\begin{lem}\label{soluble-maximal} Let $G$ be a superstable group of finite Lascar rank.

$(1)$ If $A$ and $B$ are  soluble subgroups such that $A \leq B$, then  $A^* \leq B^*$. 

$(2)$  If $(A_i)_{i \in \lambda}$ is a  totally  ordered (for inclusion) family of  $*$-connected  soluble subgroups of $G$, then there exists a soluble $*$-connected group $B$ such that $\<A_i|i \in \lambda\>\leq B$ and $B$ is  a finite intersection of $\varphi_k$-definbale soluble groups ($0 \leq k \leq 2n$),  and  $B$  is in particular definable. 

$(3)$ If $A$ is a soluble subgroup, there exists a  $*$-connected   soluble group $B$, which is maximal among $*$-connected soluble subgroups,  such that $A^* \leq B$,  and  $B$ is  a finite intersection of $\varphi_k$-definbale soluble groups ($0 \leq k \leq 2n$). 

$(4)$  For any $g \in G$, of infinite order,   there exists a  $*$-connected   soluble group $B$, which is maximal among $*$-connected soluble subgroups,  such that for some $m \in \mathbb N$, $m \neq 0$,  $g^m \in B$, 
and  $B$ is  a finite intersection of $\varphi_k$-definbale soluble groups ($0 \leq k \leq 2n$). \end{lem}

\proof

$\;$

$(1)$  Clearly, $B^* \cap A$ has a finite index in $A$, and thus $A^* \leq B^*$. 

$(2)$  By Corollary  \ref{cor-locally-soluble},   $\<A_i|i \in \lambda\>$ is soluble and hence, by Lemma \ref{prop2}, it lies in a $\varphi_k$-definable soluble subgroup $B$ ($0\leq k \leq n$).  Since $A_i$ is $*$-connected, we find by (1), $A_i \leq B^*$ and thus $\<A_i|i \in \lambda\> \leq B^*$.   Since $B$ is maximal,  it is $\phi_{i}$-definable for some $i$, and thus $B^*$ is   a finite intersection of $\varphi_k$-definbale soluble groups ($0 \leq k \leq 2n$).

$(3)$  By (2), using Zorn's lemma, we find that any $*$-connected soluble subgroup lies in a  $*$-connected   soluble group $B$, which is maximal among $*$-connected soluble subgroups, and  which is  a finite intersection of $\varphi_k$-definbale soluble groups ($0 \leq k \leq 2n$).

$(4)$  Apply (3) to the subgroup generated by $g$.  \qed

\begin{defn} Let $G$ be a superstable group of finite Lascar rank, acting nontrivially and without inversions on a $\Lambda$-tree. We define an \textit{$h$-subgroup}  to be a maximal $*$-connected soluble subgroup containing a hyperbolic element. 

Recall that if $g$ is a hyperbolic element, then $g^n$ is also hyperbolic,  for $n \neq 0$. Therefore, by Lemma \ref{soluble-maximal} (4), for any   hyperbolic element $g$, there exists an $h$-subgroup $B$ such that $g^n \in B$ for some $n \neq 0$.  We notice also that an $h$-subgroup is equationally-definable. Since  a conjugate of a hyperbolic element is hyperbolic, and since in the definition of $\varphi_k(x)$ we have taken conjugates,  it follows from the definition that if $B$ is an $h$-subgroup then $B^g$ is also an $h$-subgroup. 
\end{defn}

\begin{lem}\label{lem-h-subgroups} If $G$ is a superstable group of finite Lascar rank, then the family $\mathcal B$ of $h$-subgroups is a ferment family.  
\end{lem}

\proof The fact that $\mathcal B$ is closed by taking conjugates follows from the precedent discussion. Hence we show that if $B \in \mathcal B$, then $N_G(B)/B$ is finite. 

Suppose towards a contradiction that  $N_G(B)/B$ is infinite.  By \cite[VI, Corollary 1.3]{Berline-Lascar}, $N_G(B)/B$ has an infinite abelian subgroup. Hence there exists a definable soluble subgroup $K$ such that $B \leq K$,   with $K/B$  infinite. Since  $B \leq K^*$ and  $B$ has infinite index in $K^*$  which  contains a hyeprbolic element, we conclude that $B$ is not maximal; a contradiction. \qed

\begin{prop} \label{irre-borel-family} Let $G$ be an $\alpha$-connected superstable group of finite Lascar rank with an irreducible action on a $\Lambda$-tree, where $\Lambda=\mathbb Z$ or $\Lambda=\mathbb R$.  Suppose that  if $H \leq G$ is a definable subgroup such that $U(H)<U(G)$,  and $H$ contains a hyperbolic element, then $H$ is nilpotent-by-finite.  Then  every $h$-subgroup is nilpotent,  the family $\mathcal B$ of $h$-subgroups is a Borel family and there exists $m \in \mathbb N$ such that for any hyperbolic element $g$ there exists $1\leq n \leq m$ such that $g^n \in B$ for some $B \in \mathcal B$. If  $G=G_1*_AG_2$ and the considered action is the action of $G$ on the Bass-Serre tree associated to the precedent splitting, then    $A$ has biindex $2$ in both $G_1$ and $G_2$.
\end{prop}

\proof  The proposition is a consequence of claims below. 

\bigskip
\noindent
\textit{Claim 1.  $G$ is without soluble normal subgroup containing a hyperbolic element. }

\proof 

Suppose towards a contradiction that $G$ has a soluble normal subgroup $K$ containing a hyperbolic element.  Since  the action of $G$  is irreducible, by \cite[Lemma 2.9(2)]{chiswell}, as $K$ contains a hyperbolic element, the action of $K$ is irreducible; a contradiction as an action of  a soluble  group is reducible. \qed

\bigskip
\noindent
\textit{Claim 2. An   $h$-subgroup  is  nilpotent. }

\bigskip
\noindent
\textit{Proof.} Let $B$ be an  $h$-subgroup. By Claim 1,  $U(B)<U(G)$.  Since $B$  contains a hyperbolic element and definable,    $B$ is nilpotent-by-finite. Let $N$ be a normal nilpotent subgroup of $B$ of finite index.  Then $N$ is contained in a maximal nilpotent subgroup $\varphi_k$-definable $M$ for some $n \leq k \leq 2n$. Hence $M \cap B$ has a finite index in $B$ and therefore, since $B$ is $*$-connected,  $B\leq M$. Hence $B$ is nilpotent.  \qed

\bigskip
\noindent
\textit{Claim 3. If   $B_1$ and $B_2$   are two distinct $h$-subgroup then    $B_1 \cap B_2$  does not contain a hyperbolic element. }

\bigskip
\noindent
\textit{Proof.}

Suppose towards a contradiction that $B_1  \cap B_2$ contains a hyperbolic element and choose $B_1$ and $B_2$ such that $U(B_1  \cap B_2)$ is maximal.  

Let $H_i=  N_{B_i}(B_1 \cap B_2)$, for $i=1,2$.  Since $B_1$ and $B_2$ are $*$-connected and maximal,  $B_1 \cap B_2$ has an infinite index in  both $B_1$ and $B_2$.

Therefore, by Lemma \ref{lem-nilpotent-subgroups},$$
U(H_i)>U(B_1 \cap B_2), \;\hbox{ for } i=1,2.  \leqno (1)
$$

Let $H=N_G(B_1 \cap B_2)$. Since $B_1 \cap B_2$ contains a hyperbolic element, by Claim 1, $H$ is a proper definable subgroup of $G$. Since it contains   a hyperbolic element, it is nilpotent-by-finite.  

Let $F(H)$ be the Fitting subgroup of $H$.  Then $B_1 \cap B_2 \leq F(H)$, 
and $F(H)$ has a finite index in $H$.

Since $H$ contains some hyperbolic element $g$,  $g^n \in F(H)$ for some $n \neq 0$; and thus $F(H)$ contains  a hyperbolic element.  Therefore, by Lemma \ref{soluble-maximal} (3),  there is an $h$-subgroup $B$ such that $F(H) \cap  B$ has a finite index in $F(H)$.  

Therefore,  using (1) 
$$
U(B \cap B_1) \geq U(F(H) \cap B_1) = U(H \cap B_1) \geq U(H_1 \cap B_1)=U( H_1) > U(B_1 \cap B_2),  
$$
and 
$$
U(B \cap B_2) \geq U(F(H) \cap B_2) = U(H \cap B_2) \geq U(H_2 \cap B_2)=U( H_2) > U(B_1 \cap B_2), $$
we find 
$$
U(B \cap B_1) > U(B_1 \cap B_2),  \; U(B \cap B_2)> U(B_1 \cap B_2).
$$
But since  $B_1$ and $B_2$ are distinct,  $B \neq B_1$ or $B \neq B_2$; a contradiction with the choice of $B_1$ and $B_2$. This ends the proof of the claim. \qed

\bigskip
\noindent
\textit{Claim 4.  $\mathcal B$ is a Borel family and   there exists $m \in \mathbb N$ such that for any hyperbolic element $g$, there exists $1 \leq n \leq m$,  such that $g^n$ is in some  $B \in \mathcal B$. }

\bigskip
\noindent
\textit{Proof.} The fact that $\mathcal B$ is a Borel family is a consequence of Proposition \ref{prop-conjuguaison-borel}, Lemma \ref{lem-h-subgroups}, Lemma \ref{lem-hyper3},  Claim 3 and of the fact that an action of a nilpotent group on a real tree is abelian.

Let $g$ be a hyperbolic element. Then there exists $n \geq 1$ and an $h$-subgroup $B$ such that $g^n \in B$.  Hence $g^n \in B \cap B^g$ and thus, by Claim 3, $g \in N_G(B)$.  Since $N_G(B^h)=N_G(B)^h$,  from the fact that $h$-subgroups are conjugate, we conclude that there exists $m \geq 1$ such that for any $B \in \mathcal B$, $|N_G(B)/B| \leq m$.   \qed

\smallskip

We now suppose that $G=G_1*_AG_2$ and we consider the action of $G$ on the Bass-Serre tree associated to the precedent splitting.

\bigskip
\noindent
\textit{Claim 5.  There exist $g_1 \in G_1 \setminus A$, $g_2 \in G_2 \setminus A$ such that $G_1=A \cup Ag_1A$ and $G_2=A \cup Ag_2A$.  Thus $A$ has biindex $2$ in both $G_1$ and $G_2$. }

\bigskip
\noindent
\textit{Proof.}

Let $B$ be an $h$-subgroup. By Lemma \ref{thm3},  $B=\<C,s|C^s=C\>$, where $C \leq Stab(x)$ and $s$ is hyperbolic. After conjuguation, if necessary, we may assume that $s$ is  c.r. 

Let us denote by $|.|$ the natural length function of the free product $G=G_1*_AG_2$. 
Since $C$ is elliptic, $C \leq G_1^g$ or $C \leq G_2^g$ for some $g \in G$. Hence there exists $p$ such that for any $c \in C$, $ |c| \leq p$. 

Let $n=|s|+p$ and $ m \geq 4n$. Let $g \in G_1\setminus A$, $h \in G_2 \setminus A$. Let 
$
\hat{g}=(gh)^m. 
$
Clearly $\hat g$ is a c.r. element and it is hyperbolic. 

By Claim 5,   there is  an element $t \in G$, and $q \in \mathbb N$, $q \neq 0$,  such that $(\hat{g}^q) ^t \in B$. Hence,  since $(\hat{g}^q)^t$ is hyperbolic, we have $(\hat{g}^q)^t=as^{q'}$ for some $q' \in \mathbb Z$, with $q' \neq 0$, and $a \in C$.

We can write $(\hat{g}^q)^t$ as $d^{t'}$, where $d$ is a c.r conjugate of $\hat{g}^q=(gh)^{mq}$ and such that there is no cancellation in the products $t'^{-1}d$ and $dt'$. We notice that $|d|=|(gh)^{mq}|$. 
Since $s^{q'}= a^{-1}d^{t'}$ and $$mq \geq m \geq 4n >4|a| \hbox{ and }mq \geq m \geq 4n >4|s|,$$    we find that  $s^{\pm 1}= \alpha s'\beta$, where $\alpha, \beta \in A$ and $s'$ is the product of a subword of the normal form of $d$.

Write $s=s_1 \cdots s_q$ in normal form. We suppose without loss of generality that $s_1 \in G_1 \setminus A$.  By  properties of normal forms, we find that one of the following cases holds:
$$
s_1=\alpha_1 g \beta_1, \quad  s_2=\alpha_2 h \beta_2, \hbox{ where } \alpha_i, \beta_i \in A, \leqno (1)   
$$
$$
s_1^{-1}=\alpha_1 g \beta_1, \quad  s_2^{-1}=\alpha_2 h \beta_2, \hbox{ where } \alpha_i, \beta_i \in A, \leqno (2)
$$
and if $q \geq 3$, 
$$
s_3=\alpha_3 g \beta_3, \hbox{ where } \alpha_3, \beta_3 \in A, \hbox{ whenever }(1) \hbox{ holds },   \leqno (3) 
$$
$$
s_3^{-1}=\alpha_3 g \beta_3, \hbox{ where } \alpha_3, \beta_3 \in A, \hbox{ whenever } (2) \hbox{ holds}.  \leqno (4) 
$$

It follows in particular that, if $q \geq 3$, $$s_1, s_3 \in AgA \hbox{ or } s_1, s_3 \in Ag^{-1}A. \leqno (5) $$  

\smallskip
If we take others elements $g' \in G_1 \setminus A$, $h' \in G_2 \setminus A$, we obtain the same conclusion :  $s_1^{\pm1}=\alpha_1' g' \beta_1'$ and $s_2^{\pm 1}=\alpha_2' h' \beta_2'$, where $\alpha_i', \beta_i' \in A$.

Therefore, by (1) and (2),  for any $g' \in G_1$, $g' \in A \cup AgA \cup Ag^{-1}A$.

Now we treat the case $q=2$.  In that case we replace  $\hat g$ by $(ghg^{-1}h)^m$.  The same method yields one of the following cases
$$
s_1= \alpha_1 g \beta_1,  \quad  s_2=\alpha_2 h \beta_2, \quad s_1= \alpha_3 g^{-1} \beta_3, \leqno (6) 
$$
$$
s_1^{-1}= \alpha_1 g \beta_1,  \quad  s_2^{-1}=\alpha_2 h \beta_2, \quad s_1^{-1}= \alpha_3 g^{-1} \beta_3, \leqno (7) 
$$
which imply $G_1= A \cup AgA$.

We treat the case $q \geq 3$.  By (5) either $s_1, s_3 \in AgA \hbox{ or } s_1, s_3 \in Ag^{-1}A$. We treat only the case $s_1, s_3 \in AgA$, the other case can be treated similarly.

Then, by  replacing $\hat g$ by $( ghg^{-1}h)^m$, in the above argument, we obtain one of the following cases:
$$
s_1= \alpha_1 g \beta_1,  \quad  s_2=\alpha_2 h \beta_2, \quad s_3= \alpha_3 g^{-1} \beta_3, \leqno (8) 
$$
$$
s_1= \alpha_1 g^{-1} \beta_1,  \quad  s_2=\alpha_2 h \beta_2, \quad s_3= \alpha_3 g \beta_3, \leqno (9) 
$$
$$
s_3^{-1}= \alpha_1 g \beta_1,  \quad  s_2^{-1}=\alpha_2 h \beta_2, \quad s_1^{-1}= \alpha_3 g^{-1} \beta_3, \leqno (10)
$$ 
$$
s_3^{-1}= \alpha_1 g^{-1} \beta_1,  \quad  s_2^{-1}=\alpha_2 h \beta_2, \quad s_1^{-1}= \alpha_3 g \beta_3, \leqno (11)
$$

which imply $AgA=Ag^{-1}A$. 

By the symmetry of the problem, we conclude also that $G_2=A \cup AhA$ for some $h \in G_2$. This ends the proof of the claim as well as that of the theorem. 
\qed

\section{Proofs of the main results.}

\smallskip
\noindent
\textbf{Proof of Theorem \ref {thm1}.}

Suppose that $G$ is   $\omega$-stable and acts nontrivially on a simplicial  tree $T$. Replacing $T$ by its barycentric subdivision, we may assume that $G$ acts without inversions. Let $g$ be a hyperbolic element and   $C=Z(C_G(g))$. Then $g \in C$ and the action of $C$ on $T$, by restrection,  is nontrivial and  without inversions. Hence by Lemma \ref{lem-locally-nilpotent},   $C=\<A,s|A^s=A\>$ for some subgroup $A$.

Since $C$ is abelian and $\omega$-stable, by  Theorem \ref{macintyre-abelian}, $C=B \oplus D$, where $B$ has a finite exponent and $D$ is divisible.  There are  $b \in B$ and $d \in D$  such that $s=bd$. Since $s$ has an infinite order, $d \neq 1$.  Since $D$ is divisible, we get $d=d'^2$ for some $d' \in D$. Thus $1=|s|=|bd'^2|=|d'^2|>|d'| \neq 0$;  where here $|x|$ denotes the length of normal forms of $x$ relatively to the HNN-decomposition of $C$. A contradiction which  completes the proof.  \qed

\smallskip
\noindent
\textbf{Proof of Theorem \ref {thm2}.}

If $G$ is $\alpha$-connected and $\Lambda=\mathbb Z$, then its action is irreducible by Corollary \ref{cor-irreducible}. So we show the theorem for superstable groups with  irreducible actions.  

Let $H$ be a subnormal subgoup, such that the action of $H$ by restrection is irreducible and such that $U(H)$ is minimal.  Hence, if $K$ is definable and subnormal  with $U(K)<U(H)$, then the action of $K$ by restrection is reducible.  By Proposition \ref{irreducibles-actions}, there exists $H_1 \lhd H_2 \lhd H$ such that $H_2/H_1$ is simple and acts nontrivially on a $\Lambda$-tree, which is the desired conclusion.  \qed

\smallskip
\noindent
\textbf{Proof of Corollary \ref{cor2} and \ref{cor3}.}

Corollary \ref{cor2} is an ammediate consequence of Theorem \ref{thm1} and Theorem \ref{thm2}. Corollary \ref{cor3} is a consequence of Theorem \ref{thm1}, Theorem \ref{thm2} and Lemma \ref{lemme-dihedral} where we consider the action of $G$ on the Bass-Serre tree associated to the indicated splitting. \qed

\smallskip
\noindent
\textbf{Proof of Theorem \ref {thm3}.}

Let $G$ be a superstable group of  finite Lascar rank, having a nontrivial action  on some $\Lambda$-tree where $\Lambda=\mathbb Z$ or $\Lambda=\mathbb R$.  Suppose that,  if $H$ is a  definable  subgroup such that $U(H)<U(G)$, and having a nontrivial action on some $\Lambda$-tree, then $H$ is  nilpotent-by-finite.

We are going to treat the following cases: 

\smallskip
(1) $G^\circ$ is of infinite index in $G$. 

(2) The action of $G$ is reducible and $G^\circ$ has a finite index in $G$. 

(3) The action of $G$ is irreducible and $G^\circ$ has a finite index in $G$.

\bigskip
\noindent  \textit{\textbf{Case 1}.  $G^\circ$ is of infinite index in $G$.  } 

\smallskip 
By Lemma \ref{thm-connected},  $G$ has a definable connected normal subgroup $H \leq G^\circ$   such that  $G/H$ is nilpotent-by-finite. If $H$ acts nontrivially  on some $\Lambda$-tree, then $H$ is nilpotent-by-finite, and since it is connected, it is nilpotent.  Therefore $G$ is soluble-by-finite. By  \cite[Theorem 1.1.10]{wagner-stable}, $G$ has a definable soluble normal subgroup $K$ of finite index.  By taking $H_1=1$ and $H_2=K$, we get conclusion (1) of the theorem.

We now assume that any action  of $H$ on a $\Lambda$-tree is trivial.  We suppose also that $H$ is infinite, as otherwise $H=1$ and thus  $G$ is nilpotent-by-finite and the conclusion follows as before. 

By Lemma \ref{quotient1-actions}, $G/H$ acts nontrivially  on some $\Lambda$-tree.  Since $G/H$ is nilpotent-by-finite, by  \cite[Theorem 1.1.10]{wagner-stable}, $G/H$ has a definable nilpotent normal subgroup $K$ of finite index. By Lemma \ref{quotient1-actions}, $K/H$ acts   nontrivially  on some $\Lambda$-tree. By taking $H_1=H$ and $H_2$ to be the preimage of $K$, we obtain  conclusion (1) of the theorem. \qed

\bigskip
\noindent  \textit{\textbf{Case 2}.  The action of $G$ is reducible and $G^\circ$ has a finite index in $G$.  } 

\smallskip
It follows that $G^\circ$ is definable and has a reducible action on a $\Lambda$-tree. By Proposition \ref{reducible-actions-on-reel-trees}, the action of $G^\circ$ is abelian and there is a definable  normal subgroup $K$ such that $G/K$ is metabelian and has a nontrivial action on a $\Lambda$-tree.  If $K$ has a nontrivial action on a $\Lambda$-tree then $K$ is nilpotent-by-finite and the conclusion follows as before. If any action of $K$ on a $\Lambda$-tree is trivial, then by taking $H_0=K$ and $H_1=G^\circ$   we get the desired conclusion. \qed

\bigskip
\noindent  \textit{\textbf{Case 3}.   The action of $G$ is irreducible and $G^\circ$ has a finite index in $G$.   } 

\smallskip
It follows as above,  that $G^\circ$ is definable and has an  irreducible action on a $\Lambda$-tree.  Hence by Proposition  \ref{irreducibles-actions} we get  $K \lhd  G^\circ$ such that $K \subseteq N$, and $G^\circ/K$ is simple  with an irreducible action   on a $\Lambda$-tree. Now $G^\circ/K$ satisfies all the assumptions of Proposition \ref{irre-borel-family}, and hence we find conclusion (2) of the theorem by taking $H_1=K$ and $H_2=G^\circ$. \qed

\smallskip 

\section{Remarks}
\textbf{(1)}   Among groups acting nontrivially and without inversions on trees, we have free products of two nontrivial groups. So it is natural to see first if such groups can be superstable. The next proposition was proved by Poizat \cite{poizat-type-generic} using a beatiful technic of generic types.  We provide a  proof which is accessible to non-logicians (the main idea is that a nontrivial free product is not simple).

\begin{prop} If a nontrivial free product $G_1*G_2$ is superstable then  $G_1=G_2=\mathbb Z_2$. 
\end{prop}

\proof 

Let $1=H_0 \lhd H_1 \lhd \dots \lhd H_n \lhd G$ be the serie given  by Theorem \ref{thm-Baudisch}. Hence,  $G$ has  an infinite subnormal  subgroup $H_1$ such that either $H_1$ is abelian or $H_1/Z(H_1)$ is simple, where  $Z(H_1)$ is finite. 

We treat the first case.  By \cite[Corollary 1.2.12]{wagner-stable}, $H_1$ lies in a normal nilpotent subgroup. By   Lemma \ref{normal-small}, the trivial group is of index 2 in both $G_1$ and $G_2$ and this gives the desired conclusion. 

We treat now the second case and we show that we get a contradiction. We show by induction on $1\leq i \leq n$, that $H_i$ is in some factor. By the Kurosh subgroup theorem \cite[Theorem 1.10, Ch IV]{LyndonSchupp77},  $H_1$  is a free product $H_1=C_1 *\cdots *C_n *F$ where $F$ is a free group and each $C_i$ is the intersection of $H_1$ with some conjugate of a factor.  Therefore either $Z(H_1)=1$ or $Z(H_1)$ is infinite cyclic. Since the later case is excluded, we find that $Z(H_1)=1$. Hence $H_1$ is simple, and therefore it is in some conjugate of a factor.  

Using similar arguments, by induction on $1\leq i \leq n$, each $H_i$ is  in some conjugate of a factor.  Hence $G$ has a conjugate of a factor of finite index, which is clearly a contradiction. \qed

\smallskip
\textbf{(2)}  We give a short proof of Theorem  \ref{thm-Baudisch} by using Theorem \ref{theo-alter}.  

\proof  The proof is by induction on $U(G)$.  Let $G$ be superstable of Lascar rank $U(G)= \omega^\alpha n+\beta$ where
$\alpha$ and $\beta$ are ordinals, $n \in \mathbb N^*$ and $\beta
< \omega^\alpha$, and suppose that the theorem holds for all superstable groups $H$  such that $U(H)<U(G)$.

Suppose that $G$ has a definable normal subgroup $K$  such that $ \omega^\alpha \leq U(K)<U(G)$.  By induction, $K$ and $G/K$ satisfies the conclusions of the theorem. 

Since, in our context, a finite-by-(definable abelian) group is (definable abelian)-by-finite and a   finite-by-(definable simple) group is centre-by-(definable simple), the conclusion for $G$ follows by induction. 

So we suppose that if $H$ is a definable normal subgroup of $G$, then  either $U(H)<\omega^\alpha$ or $H$ has a finite index in $G$. 

By Theorem \ref{theo-alter}, $G$ has either a normal sbgroup $H$ such that $U(H) \geq \omega^\alpha$ and $H$ is either $\alpha$-connected or nilpotent. 

If $H$ is nilpotent then $G$ is nilpotent-by-finite, and by \cite[Corollary 1.1.11]{wagner-stable}, $G$ satisfies the desired conclusions. 

If $H$ is $\alpha$-connected, then we replace $G$ by $H$ and thus we may assume that any definable normal subgroup of $G$ satisfies $U(H)<\omega^\alpha$. Therefore, by \cite[VI, Proposition 2.6]{Berline-Lascar}, $G/Z(G)$ is simple or abelian and thus we get the desired conclusion. \qed

\smallskip
\textbf{(3)}  Here a proposition analogous to Corollary \ref{quotient-par-la-composante-connexe}. 

\begin{prop} \label{prop-component} If $G$ is a superstable group then $G/G^{\circ}$ is soluble-by-finite. 
\end{prop}

We need the following lemma. 

\begin{lem} \label{lem-resid}If $G$ is a residually finite group and $H$ is a soluble normal subgroup of $G$, then there exists a soluble subgroup $K$ such that $H \lhd K \lhd G$ and $G/K$ is residually finite. 
\end{lem}

\proof  Let us say that a subgroup $H$ of a group $G$ is $\bigwedge$-e-definable if $H$ can be defined as an  intersection of equations; that is  if
there exist  words $(w_i(x)| i \in I)$, with parameters from
$G$,  such that $H=\bigcap_{i \in I}\{g \in G\mid G \models
w_i(g)=1\}$.

We   first prove the following general property: \textit{if $H \lhd G$ is a normal soluble subgroup, then $H \lhd K \lhd G$, with $K$ is soluble and $\bigwedge$-e-definable}. 

The proof is by induction on the derived length of $H$. If $H$ is of derived length $0$, the result is clear. 

Suppose that $H$ is of derived lenght $n+1$. 

Then $H^{(n)}$ is abelian
and normal. Let $N$ be the center of $C_G(H^{(n)})$ in $G$. Then
$H^{(n)} \leq N \trianglelefteq G$. 

Since $N$ can be written as an
intersection  of centralizers of elements of $G$, $N$ is   $\bigwedge$-e-definable. 

Let $\bar G=G/N$ and  $\pi : G \rightarrow \bar G$ be the canonical morphism. Then $\pi(H)$
is a normal and solvable  subgroup of derived length $n$ in $\bar G$.
Thus, by induction, $\pi(H) \leq K_1$ and $K_1$ is a solvable normal
and $\bigwedge$-e-definable subgroup of $\bar G$. Then $H \leq \pi^{-1}(K_1) \lhd G$, and $\pi^{-1}(K_1)$ is soluble. 

Since $K$ is  $\bigwedge$-e-definable in $G$ and $K_1$ is $\bigwedge$-e-definable subgroup in $\bar G$, it follows that $\pi^{-1}(K_1)$ is  $\bigwedge$-e-definable in $G$. This ends the proof of the enounced general property. 

Now we show that if $G$ is a residually finite group and $K \lhd G$ is  $\bigwedge$-e-definable in $G$, then $G/K$ is residually finite.  The proof proceeds as in \cite[Lemma 2.2]{ould-super-resid}. 

 Let $K$ be definable by $\bigwedge_{i \in I}w_i(\bar
a_i, x)=1$ and   $\pi : G \rightarrow G/K$ be the canonical
morphism.

Let $g \in G$ such that $\pi(g) \neq 1$; that is $g \not \in K$.
Then there is some $p \in I$ such that $w_p(\bar a_p, g)
\neq 1$. Since $G$ is residually finite, there exists a surjective
morphism $\phi : G \rightarrow L$, where $L$ is finite,  such that
 $\phi(w_p(\bar a_p, g)) \neq 1$.

We claim that $\phi(g) \not \in \phi(K)$. If $\phi(g) \in
\phi(K)$,  then there exists an element $g' \in K$ such that
$\phi(g)=\phi(g')$. Since $g' \in K$ we get $\phi(w_p(\bar a_p,
g'))=w_p(\phi(\bar a_p), \phi(g'))=1$  and thus $\phi(w_p(\bar
a_p, g))=1$; a contradiction.

Since $\phi$ is surjective, $\phi(K)$ is a normal subgroup of $L$.
Let $\pi' : L \rightarrow H=L/\phi(K)$ be the canonical morphism.
We define $f : G/K \rightarrow H$ by $f(\pi(x))=\pi'(\phi(x))$.
Then $f$ is a morphism.   Now if $f(\pi(g))=1$, then
$\pi'(\phi(g))=1$ and thus $\phi(g) \in \phi(K)$; a contradiction.
Thus $G/K$ is residually finite as desired. \qed

\smallskip
\noindent 
\textbf{Proof of Proposition \ref{prop-component}}.  

We  first show that $\bar G=G/G^\circ$ is residually finite.  Let $\pi : G \rightarrow \bar G$ be the natural homomorphism and $g \in G$, $\pi(g) \neq 1$.  Since $g \not \in G^{\circ}$ and $G^{\circ}$  is an intersection of definable normal subgroups of finite index, there exists a finite index normal definable subgroup  $H$ such that $g \not \in H$. Since $G^\circ \leq H$, we get  $\pi(g) \not \in \pi(H)$, and thus we have the desired conclusion.

Let  $1=H_0 \lhd H_1 \lhd \dots \lhd H_n \lhd G$ be the serie given by Theorem \ref{thm-Baudisch}.   Hence 
$\pi(H_0)  \lhd \pi(H_1) \lhd \dots \lhd \pi(H_n) \lhd \bar G$. Let $p\leq n$ be the greatest integer such that   $\pi(H_p)$ is soluble. We are going to show that $p=n$. Suppose towards a contradiction that $p<n$. 

Suppose first that $H_{p+1}/H_{p}$ is abelian.  Then $\pi(H_{p+1})/\pi(H_p)$ is abelian, hence $\pi(H_{p+1})$ is soluble; a contradiction with the choice of $p$. 

Thus $H_{p+1}/H_{p}$ is infinite simple modulo a finite centre. Since $\pi(H_{p+1})/\pi(H_p)$ is nontrivial,  there exists a soluble subgroup $K \lhd \pi(H_{p+1})$ such that $\pi(H_{p+1})/K$ is infinite simple.  By Lemma \ref{lem-resid}, and since $K$ is a maximal normal subgroup in $\pi(H_{p+1})$, $\pi(H_{p+1})/K$ is residually finite. A contradiction with the fact that  $\pi(H_{p+1})/K$ is infinite simple. Therefore $p=n$ and thus $\bar G$ is soluble-by-finite. \qed

Note that we have proved more, namely that a residually finite quotient of $G$ is soluble-by-finite. 

\smallskip
\textbf{(4)} Recall that a subgroup $H$ of a group $G$ is {\em
conjugately separated} in $G$, or {\em malnormal} in $G$, if $H
\cap H^x =1$ for every $x \in{G\setminus H}$; and a {\em CSA-group}
 is a group in which every
maximal abelian subgroup is malnormal. In dealing with actions on trees, the CSA-property has several intersting consequences. In the case of CSA-groups, with the superstable assumption, we have the following.

\begin{prop} \label{conjugacy-centralizers} Let $G$ be an $\alpha$-connected superstable CSA-group. If $C$ is a  maximal abelian subgroup then $G =\bigcup_{g \in G}C^g$. If $G$ has a nontrivial action on a simplcial tree then $G$ splits as $G=G_1*_AG_2$ with $A$ has biindex $2$ in both $G_1$ and $G_2$. 
\end{prop}

\proof 

We may assume that $G$ is nonabelian, as otherwise the result is clear. By \cite[VI, Corollary 1.3]{Berline-Lascar}, $G$ has an infinite abelian subgroup.  Hence there exists an element $g \in G$ such that  $C=C_G(g)$ is  infinite.  

We claim that  $C$ is generous.  Since $G$  has a monomial Lascar rank,  we can apply Corollary  \ref{cor-generic}. Since $G$ is a CSA-group and nonabelian, $X(C)=1$ and $C$ is selfnormalizing. Since $C$ is infinite, $X(C)$ is not generic in $C$ and hence $C$ is generous as claimed.

Let $h \in G$ and suppose that $C_G(h)$ is finite. Then $h^G$ is generic, and thus it meets $\bigcup_{g \in G}C^g$. Thus $h$ lies in an infinite abelian subgroup, a contradiction. 

Hence  any $h \in G$ has an infinite centralizer. It follows that for any $h \in G$, $C_G(h)$ is generous. Since $G$ is connected, we conclude as in the proof of Proposition \ref{prop-conjuguaison-borel}, that   any two nontrivial centralizers are conjugate.  Hence $G=\bigcup_{g \in G}C^g$ for any maximal abelian subgroup as desired.

 If $G$ has a nontrivial action on a simplcial tree, where the action is assumed to be, without loss of generality,  without inversions,   then the action  of $G$ is irreducible by Corollary \ref{cor-irreducible}. By Theorem \ref{thm-bass-serre}, $G$ splits as $G=G_1*_AG_2$. The proof of the fact that  $A$ has biindex $2$ in both $G_1$ and $G_2$ procced in a similar way to that of Proposition \ref{irre-borel-family}(Claim 5). \qed 

\smallskip
\textbf{(5)}  How about free actions of superstable groups on real trees? We have seen that $\mathbb Q$ (also $\mathbb Z$) is superstable and acts freely on a real tree. We claim that if $G$ is a superstable group acting freely on a real tree then $G$ is abelian.  By \cite[Proposition 5.5.13]{chiswell}, $G$ is locally fully residually free. Hence $G$ is a model of the universal theory of nonabelian free groups, and since it is superstable it is abelian by \cite[Corolary 1.2]{Ould-csa-super}.

\smallskip
\textbf{(6)}  What can be said about infinite finitely generated   superstable groups? It is noticed in \cite[Proposition 3.2]{ould-super-resid} that the existence of an infinite  finitely generated   $\omega$-stable  group implies the existence of a simple $\omega$-stable finitely generated one.  In the superstable case, the situation is different because the existence of finitely generated abelian superstable groups. However, in the presence of the $\alpha $-connectedness, we have   at least the following.

\begin{prop}
If $G$ is an infinite  finitely generated $\alpha$-connected superstable group then $G$ has a definable normal subgroup $N$  such that $G/N$ is infinite simple (and of course finitely generated and not algebraic over an algebraically closed field). 
\end{prop}

\proof 

Since $G$ is finitely generated, it has a maximal normal proper subgroup $N$.   We claim that $N$ is definable. If it is not the case, then $N$ has a finite index by Lemma \ref{cor-maximal-normal}. Therefore, by Lemma \ref{lem-finite-index},  $G$ has a definable normal subgroup $K \leq N$ such that $G/K$ is abelian. Since $G$ is connected and $N$ is proper, we get that $G/K$ is infinite. Since $G/K$ is  infinite abelian and finitely generated, it is not connected; a contradiction. Hence $N$ is definable as claimed and thus  $G/N$ is simple and infinite. \qed

\bigskip
\noindent {\bf Acknowledgements }  The author wishes to express his thanks to Frank Wagner for some indications  in the proof of Proposition  \ref{propo-quotient} and for some references, and to Gabriel Sabbagh for some remarks and suggestions.  
\bibliographystyle{alpha}
\bibliography{biblio}
\end{document}